\documentclass[12pt]{article}
\usepackage{mathtools} 
\usepackage{amsmath, amssymb} 
\usepackage{amsthm}
\usepackage{enumerate}  
\usepackage{url} 

\usepackage{graphicx}

\newcommand\cx{{\mathbb C}}

\newcommand\re{{\mathbb R}}

\newcommand\cS{{\cal S}}
\DeclareMathOperator\dist{dist} 
\newcommand\De{\Delta}
\newcommand\be{\beta}
\newcommand\ga{\gamma}

\DeclarePairedDelimiter\abs{\lvert}{\rvert}%
\DeclarePairedDelimiter\norm{\lVert}{\rVert}%

\makeatletter
\let\oldabs\abs
\def\abs{\@ifstar{\oldabs}{\oldabs*}}
\let\oldnorm\norm
\def\norm{\@ifstar{\oldnorm}{\oldnorm*}}
\makeatother

\newcommand\opk[1]{\mathop{\mathrm{#1}}\nolimits}

\newcommand\sbs{\subseteq}
\newcommand\sps{\supseteq}

\newcommand\comp[1]{{\mkern2mu\overline{\mkern-2mu#1}}}
\newcommand\pmat[1]{\begin{pmatrix} #1 \end{pmatrix}}
\newcommand\seq[4]{#1_{#2},#1_{#3},\ldots,#1_{#4}}

\newtheoremstyle{plainsl}%
	{\topsep}
	{\topsep}
	{\slshape} 
	{}
	{\normalfont\bfseries}
	{.}
	{ }
	{}

\swapnumbers

{\theoremstyle{plainsl}
\newtheorem{theorem}{Theorem}[section]
\newtheorem{lemma}[theorem]{Lemma}
\newtheorem{corollary}[theorem]{Corollary}}
{\theoremstyle{remark}
}

\renewcommand\proof{\noindent\textsl{Proof. }}
\newcommand\sqr[2]{{\vbox{\hrule height.#2pt
    \hbox{\vrule width.#2pt height#1pt \kern#1pt
        \vrule width.#2pt}\hrule height.#2pt}}}
\renewcommand\qed{%
	\ifmmode\eqno\sqr53
	\else\nolinebreak\ \hfill\sqr53\medbreak\fi}

\DeclareMathOperator{\Ch}{Ch}
\DeclareMathOperator{\op}{op}
\DeclareMathOperator{\Br}{Br}

\newcommand\one{{\bf1}}

\newcommand\sym[1]{\opk{Sym}(#1)}

\title{An Introduction to the Moebius Function}
\author{C. D. Godsil}

\begin{document}
\maketitle

\begin{abstract}
This is an introduction to the M\"obius function of a poset. The chief novelty is in the exposition. 
We show how order-preserving maps from one poset to another can be used to relate their 
M\"obius functions. We derive the basic results on the M\"obius function, applying them in particular 
to geometric lattices.
\end{abstract}

\tableofcontents

\section{Posets and Matrices}
Our first four sections provide a fairly standard approach to the M\"obius function of a poset. It is based in part on the treatment in Chapter~2 of Lov\'asz \cite{Lovasz1979}.

Let $P$ be a poset with elements $\seq{p}{1}{2}{n}$. (Unless we explicitly say otherwise, all posets 
we consider are finite. So $n$ is an honest-to-God Kroneckerian integer.) An $n\times n$ matrix $B$ 
is \textsl{compatible} with $P$ if $(B)_{ij}$ is zero unless $p_i\le p_j$. It is immediate that set 
of all $n\times n$ matrices compatible with $P$ is closed under addition, and it is not hard to 
show that it is also closed under multiplication. Thus it is an algebra over $\cx$, often called 
the \textsl{incidence algebra} of $P$. We note that it contains the idenity matrix, as well as 
the \textsl{Zeta matrix} $Z_p$, which has $ij$-entry equal to one if and only if $p_i\le p_j$. 
Any matrix compatible with $P$ can be regarded as a function on $P\times P$. This function is 
non-zero only on the ordered pairs $(x,y)$ where $x\le y$, and so we may even view our function as 
a function on the \textsl{intervals} of $P$.

A simple induction argument shows that, by relabelling the elements of $P$ if needed, we may assume 
that $i\le j$ whenever $p_i\le p_j$. Then the matrices compatible with $B$ are all upper triangular, 
and so such a matrix is invertible if and only if its diagonal entries are all non-zero. We have 
the following intersting result.

\begin{lemma}
	Let $P$ be a poset. If $B$ is compatible with $P$ and invertible, then $B^{-1}$ is compatible 
	with $P$.
\end{lemma}

\proof
Let $\phi(x)$ be the characteristic polynomial of $B$. If $B$ is invertible then $\phi(0)\ne 0$ and so
\[
	\phi(x) = x\psi(x) +c
\]
for some polynomial $\psi$ and non-zero constant $c$. By the Cayley-Hamilton theorem
\[
	0 =\phi(B) = B\psi(B) + cI,
\]
whence $c^{-1}\psi(B)= B^{-1}$.\qed

Since the zeta matrix $Z_p$ has all its diagonal entries equal to one, it is invertible. 
By the lemma, $(Z_p)^{-1}$ is compatible with $P$. The corresponding function on $P\times P$ 
is the \textsl{M\"obius function} of $P$, and is denoted by $\mu_P$. We also use $M_P$ to 
denote $ZP^{-1}$, thus $(M_P)_{a,b} = \mu_P(a,b)$.

We can determine $\mu_P$ by inverting the triangular matrix $Z_P$; this represents no intellectual
challenge and can be carries out in polynomial time. However, for many interesting posets, properties of
the M\"obius function can be read off from properties of the posets. The values taken by the M\"obius
function may have combinatorial significance.

\section{M\"obius Inversion}
Our first result is known as the principle of M\"obius inversion.

\begin{theorem}
	Let $P$ be a poset and let $f$ and $g$ be functions on $P$. Then 
	\begin{enumerate}[(a)]
		\item $g(x) = \sum_{y\ge x} f(y)$ if and only if $f(z) = \sum_y \mu(z,y)g(y)$.
		\item $g(x) = \sum_{y\le x} f(y)$ if and only if $f(z) = \sum_y \mu(y,z) g(y)$.
	\end{enumerate}
\end{theorem}
\proof
We may abuse notation and view $f$ and $g$ as column vectors, with entries indexed by $P$. 
Then (a) says that 
\[
	g = Z_Pf \Leftrightarrow M_P g= f
\]
and (b) that
\[
	g = Z_P^Tf \Leftrightarrow M_P^Tg = f.
\]
Since $M_P = Z_P^{-1}$, no more need be said.\qed 

Since all diagonal entries of $Z_P$ are equal to one, it follows that the same is true for $M_P$. (One way to convince yourself of this is to recall that the diagonal entries of a triangular matrix are its eigenvalues, and that the eigenvalues of $Z_P^{-1}$ are the reciprocals of the eigenvalues of $Z_P$.) Thus $\mu_P(x,x)=1$, for any element $x$ of $P$. There is a recursive expression for the remaining values of $\mu_P$, equivalent to the back-substitution phase in Gaussian elimination.

\begin{lemma}
\label{lem:muPab}
	Let $a$ and $b$ be two elements of the poset $P$. Then
	\[
		\mu_P(a,b) = \begin{cases}
						0, \quad a\not\le b;\\
						1, \quad a=n;\\
						-\sum_{x: x\le b} \mu_P(a,x), \quad \mathrm{otherwise.}
					\end{cases}\]
\end{lemma}

\proof
If $a\not\le b$ then $(M_P)_{a,b}=0$, since $M_P$ is compatible with $P$. If $a=b$ then 
$\mu_P(a,a)=1$, as noted above. Finally, if $a<b$ then $(M_PZ_P)_{ab}=0$ and therefore
\[
 0 = \sum_{x\le b} \mu_P(a,x).
 \]
Hence $\mu(a,b)= -\sum_{x< b} \mu_P(a,x)$, as required.\qed 

The argument used in the previous proof yields another useful identity. Suppose $a$ and $b$ are elements of the poset $P$ and $a<b$. Then
\[
	\mu_P(a,b) = -\sum_{x:a<x} \mu_P(x,b).
\]

The \textsl{chain} $\mathcal{C}(n)$ is the poset with elements $0,\cdots,n$, where $i\ge j$ if $i-j$ 
is non-negative. Suppose $a$ and $b$ are elements of $\mathcal{C}(n)$ and $\mu=\mu_{\mathcal{C}(n)}$. 
If $a<b$ then $\mu(a,b)=-1$ if $b$ covers $a$, and is zero otherwise. We will use this in the 
next section to compute the M\"obius function for the poset of divisors of a given integer.

\section{Products}
The \textsl{product} of poset $P$ and $Q$ is the poset with elements $P\times Q$, where 
\[(x,y) \le _{P\times Q} (x',y')\]
if and only if 
\[x\le _P x' \text{ and } y\le_Q y'.\]
We consider two examples. Let $\mathcal{B}(n)$ be the lattices of subsets of an $n$-element set. It is
routine to verify that $\mathcal{B}(n)$ is isomorphic to the product of $n$ copies of $\mathcal{B}(1)$,
which in turn is isomorphic to $\mathcal{C}(1)$. The lattice of divisors of an integer $n$ is also
isomorphic to a product of chains. More precisely, if $p$ is prime and $n=p^r$ then the lattice of
divisors of $n$ is the chain of length $r$. If $n=p^rm$ where $m$ and $p$ are coprime then the divisor
lattice of $n$ is the product of the divisor lattice of $m$ with the chain of length $r$. Note that
$\mathcal{B}(n)$ can be regarded as the divisor lattice of a square-free integer having exactly $n$
distinct prime divisors.

Turning from examples to M\"obius functions, we have
\[Z_{P\times Q} = Z_P \otimes Z_Q,\]
whence 
\[M_{P\times Q} = M_P \otimes M_Q.\]
As an immediate consequence we have the next result.

\begin{lemma}
	If $P$ and $Q$ are posets and $(x,y)$ and $(x',y')$ are elements of $P\times Q$ then 
	\[\mu_{P\times Q}((x,y),(x',y')) = \mu_P(x,x') \mu_Q(y,y'). \tag*{\sqr53}\]
\end{lemma}

Suppose that $S$ and $T$ are subsets of some $n$-element set. Then, taken with our remarks above, the previous lemma implies that
\[
	\mu(S,T)=
		\begin{cases}
			0, \quad S\not\sbs T;\\
			(-1)^{\abs{T\backslash S}}, \quad \mathrm{otherwise}.
		\end{cases}
\]

\section{Derangements}

We now present a classical combinatorial application of the M\"obius function. A \textsl{derangement} 
is a permutation with no fixed points. We wish to compute $D_n$, the number of derangements of $n$ points.

To this end, if $S\sbs\{1,\cdots,n\}$ let $D_n(S)$ denote the number of permutations of $\{1,\cdots,n\}$ which fix each point in $S$ and no points not in $S$. (So $D_n(\emptyset)=D_n$.) Let $F_n(S)$ denote the number of permutations which fix each point in $S$. Both $F_n(S)$ and $D_n(S)$ are functions on $\mathcal{B}(n)$. We have 
\[F_n(S) = (n-\abs{S})!\]
and we will use this to compute $D_n$.

The key observation is that
\[
	F_n(S)=\sum_{T\sps S} D_n(T)
\]
whence 
\begin{align*}
	D_n(S) &= \sum_T \mu_{\mathcal{B}(n)} (S,T) F_n(T)\\
		&=\sum_{T\sps S} (-1)^{\abs{T\backslash S}} F_n(T).
\end{align*}
Assuming that $\abs{S}=k$, we may write the last sum as
\[
	\sum_{\ell=k}^n (n-\ell)! \binom{n-k}{\ell-k} (-1)^{\ell-k}
\]
and therefore 
\begin{align*}
	D_n &=\sum_{\ell=0}^n (n-\ell)!\binom{n}{\ell} (-1)^{\ell}\\
		&= n! \left(1-\frac{1}{1!} + \frac{1}{2!} - \cdots+(-1)^n\frac{1}{n!}\right)\\z
		&=\left[\frac{n!}{e}\right]
\end{align*}

\section{Posets and Chains}

A \textsl{chain} in a poset is a set of elements, any two of which are comparable. Any finite chain has
unique minimal and maximal elements. The set of all non-empty chains of the poset $P$ will be denoted by
$\Ch(P)$. This set is partially ordered by inclusion, hence is itself a poset. Our first task in this
section is to describe the relation between chains and the M\"obius function. We denote the length of
the chain $C$ by $\ell(C)$. (This is one less than the number of elements in $C$.)

For this we need another definition. If $P$ is a poset with elements $\seq{p}{1}{2}{n}$, let $Y_P$ be
the $n\times n$ matrix with $ij$-entry equal to one if $p_i<p_j$. Thus, if we have arranged things so
that $Z_P$ is triangular then $Y_P=Z_P-I$.

\begin{lemma}
\label{lem:YZM}
	Let $P$ be a poset with elements $\seq{p}{1}{2}{n}$. Then
	\begin{enumerate}[(a)]
		\item the $ij$-entry of $Y_P^m$ is the number of chains of length $m$ in $P$ with 
			least element $p_i$ and maximal element $p_j$,
		\item the $ij$-entry of $Z_P^m$ is a polynomial in $m$, and 
		\item the $ij$-entry of $M_P$ is $\sum_{C\in \Ch(P)}(-1)^{\ell(C)}$.
	\end{enumerate}
	\label{4.1}
\end{lemma}

\proof
Given that $Y_P^{m+1}=Y_P^m Y$, it is easy to prove (a) by induction on $m$. If $m>\abs{P}$ then
$Y_p^m=0$. Assuming that $Z_P=I+Y_P$, we then have \[ Z_P^m=\sum_{k=0}^m \binom{m}{k} Y_P^k. \] Since
$Y_p^m=0$ for sufficiently large $k$ and since $\binom{m}{k}$ is a polynomial in $m$ (of degree $k$), it
follows that the entries of $Z_P^m$ are polynoimals in $m$.

To prove (c) we observe that 
\[
	M_P=Z_P^{-1}=(I+Y_P)^{-1}=\sum_{k\ge 0} (-1)^k Y_P^k. \tag*{\sqr53}
\]

Lemma~\ref{lem:YZM}(c) is quite important, and it is worth recording it in a slightly different form.

\begin{lemma}[P. Hall]
\label{lem:Hall}
	If $a$ and $b$ are elements of the poset $P$ then 
	\[\mu_P(a,b)=\sum (-1)^{\ell(C)},\]
	where the sum is over all chains $C$ in $P$ with minimal element $a$ and maximal element $b$.
	\qed
\end{lemma}

We can always create a new poset from $P$ by reversing the order. The reult is a poset, $P^{\circ p}$ 
say, with the same elemtns as $P$ such that
\[
	a\le_P b \Leftrightarrow b\le_{P^{\circ p}} a.
	\]
One immediate consequence of Hall's theorem is that
\[
	\mu_P(a,b) = \mu_{P^{\circ p}} (b,a).
\]
This can often be used to deriva alternate forms of various identities, e.g., the expression 
for $\mu_P$ we gave directly following Lemma~\ref{lem:muPab} can be derived from 
Lemma~\ref{lem:muPab} in this way. 

Making use of terminology to be explained later, the $ab$-entry of $(Z_P)^m$ can be shown to be 
equal to the number of order preserving mappings from a chain of length $m$ into $P$. (The 
corresponding entry of $(Y_P)^m$ counts order preserving injections.)

\section{Simplicial Complexes}

A simplicial complex $\mathcal{S}$ on a set $\Omega$ is a set of non-empty subsets of $\Omega$ such that
if $A\in \mathcal{S}$ and $B\sbs A$ then $B\in \mathcal{S}$. (Oh well, there are two schools of thought.
Some authors choose to ake the empty set an element of any simplicial complex.) The elements of
$\mathcal{S}$ are called \textsl{faces} and the \textsl{dimension} of a face $A$ is $\abs{A}-1$. (Yes,
the empty set would have dimension $-1$.) The maximal lements of $\mathcal{S}$ are called
\textsl{simplices}. We denote the number of $k$-dimensional faces of $\mathcal{S}$ by
$f_k(\mathcal{S})$, and call it the \textsl{$k$-th level number of the complex}. The \textsl{Euler
characteristic} of $\mathcal{S}$ is defined to be \[\sum_{k\ge 0} (-1)^k f_k.\]

We consider two examples. Let $\mathcal{M}$ be a triangulation of a surface and let $\mathcal{S}$ be the
simplicial complex whose elements are the sets of vertices contained in some face of $\mathcal{M}$. (To
be more prosaic, the elements of $\mathcal{S}$ are the vertices, edges and triangles of $\mathcal{M}$.)
In this case the Euler characteristic of $\mathcal{S}$ is determined by the surface on which
$\mathcal{M}$ lies.

Our second example is $\Ch(P)$. The sunokuces are the maxial chains in $P$. If $\widehat{P}$ is obtained from $P$ by adjoining a new $0$- and $1$-element then the Euler characteristic of $\Ch(P)$ is equal to 
\[
	1 + \mu_{\widehat{P}}(0,1).
\]
To make matters more confusing, we note that every simplicial complex is a poset. We will see later that $\mathcal{S}$ and $\Ch(S)$ have the same Euler characteristic.

\section{Determinants}
The theory we describe in this section is one of the prettiest parts of the theory of the 
M\"obius function, and was developed independently by Lindstr\"om \cite{Lindstrom1969} 
and Wilf \cite{Wilf1968}.

\begin{lemma}
\label{lem:gxysum}
	Let $f$ be a function defined on the poset $P$ and set 
	\[
		g(x,y)=\sum_{z\ge x,y} f(z).
	\]
	If $G$ is the matrix with rows and columns indexed by $P$ and $xy$-entry equal to $g(x,y)$, 
	then $\det(G)=\prod_{x\in P} f(x)$.
\end{lemma}
\proof
Let $F$ be the diagonal matrix with rows and columns indexed by $P$, where $(F)_{xx}=f(x)$. 
Then $G=Z_P F Z_P^T$ and so
\[
	\det G = \det(Z_P F Z_P^T) = (\det Z_P)^2 \det F = \det F. \tag*{\sqr 53}
\]

\textbf{Exercise:} Give an expression for $f$ in terms of $g$.

If $P$ is a lattice then $g(x,y)=\sum_{z\ge x\vee y} f(z)$. Thus we may allow $g$ to be any 
function defined on $P$, with $f$ given by
\[
	f(y)=\sum_z \mu_P(y,z) g(z).
\]
Then Lemma~\ref{lem:gxysum} implies that 
\[\det G = \prod_{x\in P} \sum_{y\in P} \mu_P(x,y) g(y).\]

We will make significant use of this result later. Further applications appear in the papers 
of Lindstr\"om and Wilf \cite{Lindstrom1969,Wilf1968} and in \cite[Ch 2]{Lovasz1979}.

\section{Distance Matrices of Trees}

If $X$ is a graph on $n$ vertices, its \textsl{distance matrix} is the $n\times n$
matrix $D=D(X)$ such that
\[
	D_{u,v} := \dist_X(u,v).
\]
This a symmetric matrix with zero diagonal. In this section we determine the determinant
of the distance matrix of a tree, and derive an explicit formula for the inverse of
its distance matrix.

There is a natural construction of a partial order from a rooted tree.  Assume we have
a tree $T$ with a root vertex $x$.  If $u,v\in V(T)$, we declare that $\le v$ if
$u$ lies on the unique path from $x$ to $v$.  We define the incidence matrix $Z$
of this partial order to be the $|V(T)| \times|V(T)|$ matrix given by
\[
	(Z)_{u,v} =\begin{cases}
				1,& \textrm{if } u\le v;\\
				0,& \textrm{otherwise}.
			\end{cases}
\]
We may assume without loss that $Z$ is upper triangular.  Since $(Z)_{u,u}=1$ for
each vertex $u$ of $T$, it follows that $\det(Z)=1$ and $Z$ is invertible.

Define the square matrix $H$ by
\[
	H := \one e_1^T +e_1\one^T -2I;
\]
thus
\[
	H=\pmat{0&\one^T\\ \one&-2I}.
\]
Observe that
\[
	H\pmat{1&0^T\\ -\frac12\one&I} =\pmat{\frac12(n-1)&\one^T\\ 0&-2I}.
\]
and so
\begin{equation}
	\label{detH}
	\det(H) =\frac12(n-1)(-2)^{n-1}
\end{equation}

Graham and Lov\'asz \cite{GrahamLovasz} proved the following result.  

\begin{theorem}
	If $D$ is the the distance matrix of the tree $T$, then 
	\[
		D(T)=Z^THZ.
	\]
\end{theorem}

\proof
Assume that we have arranged for $Z$ to be upper triangular. If we denote $Z^T\one$ by $d$, then
\[
	Z^T\one e_1^T Z +Z^Te_1\one^TZ = d\one^T +\one d^T
\]
and therefore
\[
	Z^THZ =d\one^T +\one d^T -2Z^TZ.
\]
If $0$ is the root of $T$, the $ij$-entry of the right side here is equal to
\begin{align*}
	(1+\dist(0,i)) +(1+\dist(0,j)) &-2(1+\dist(0,i\wedge j))\\
		&=\dist(0,i)) +\dist(0,j) -2\dist(0,i\wedge j)\\
		&=\dist(i,j).\tag*{\sqr53}
\end{align*}

Since $\det(Z)=1$, the following famous theorem of Graham and Pollak \cite{GrahamPollak} 
is an immediate consequence.

\begin{theorem}
	\label{thm:gp}
	If $T$ is a tree on $n$ vertices, then 
	\[
		\det D(T) = (n-1)(-1)^{n-1}2^{n-2}.\qed
	\]
\end{theorem}

From the observation that $D(T)=Z^THZ$ we deduce that
\[
	D(T)^{-1} =Z^{-1}H^{-1}Z^{-T}.
\]
Since $Z$ is upper triangular and $(Z)_{i,i}=1$ for all $i$, we see that $Z^{-1}$
is an integer matrix and is also upper triangular.  In fact we can be much more precise.

\begin{lemma}
	Let $T$ be a rooted tree and let $Z$ be the incidence matrix of the
	partial order it determines.  Then
	\begin{enumerate}[(a)]
		\item
		$(Z^{-1})_{u,u} =1$ for all vertices $u$;
		\item
		If $u$ is adjacent to $v$ and lies on the unique path from the root to $v$, then
		$(Z^{-1})_{u,v}=-1$.
		\item
		Otherwise $(Z^{-1})_{u,v}=0$.
	\end{enumerate}
\end{lemma}

\proof
These are immediate conseqwuences of the fact that each interval of the poset determined 
by $T$ is a chain.\qed

We use this lemma to obtainan expression for the inversed of the distance matrix of a
tree.  The columns of $Z^{-1}$
are indexed by the vertices of $T$. The column corresponding to the root is the characteristic
vector of the root, viewed as a subset of $V(T)$.  The remaining columns are the signed
characteristic vectors of the edges of $T$.  Hence
\[
	Z^{-1} =\pmat{e_1& \tilde B}
\]
where $\tilde B$ is the incidence matrix of an orientation of $T$.  (Depending on assumptions,
it may be taken to be the orientation where each edge is directed towards the root.)

\begin{theorem}
	\label{thm:dinv}
	Let $T$ be a tree with valency matrix $\De$ and let $\be$ denote the vector
	$(2I-\De)\one$.  Then
	\[
		D^{-1} =\frac1{2n-2}\be\be^T -\frac12(\De-A).
	\]
\end{theorem}

\proof
By multiplication we can verify that
\[
	H^{-1} =\frac1{n-1}\pmat{2&\one^T\\ \one& \frac12(J-(n-1)I)}.
\]
Denote $\tilde B\one$ by $\ga$.  We have
\begin{multline*}
	\pmat{e_1&\tilde B} \pmat{2&\one^T\\ \one&\frac12(J-(n-1)I)} \pmat{e_1^T\\ \tilde B^T}\\
		=2e_1e_1^T +e_1\ga^T +\ga e_1^T +\frac12(\ga\ga^T -(n-1)\tilde B\tilde B^T).
\end{multline*}
Now
\[
	4e_1e_1^T +2e_1\ga^T +2\ga e_1^T +\ga\ga^T =(2e_1+\ga)(2e_1+\ga)^T
\]
and since $2e_1+\ga=\be$ and $\tilde B\tilde B^T=\De-A$, the result follows.\qed

\section{Order-Preserving Mappings}

Let $P$ and $Q$ be posets. A function $f$ form $P$ to $Q$ is \textsl{order-preserving} if, whenever $x$
and $y$ belong to $P$ and $x\le y$, we have $f(x)\le f(y)$. To consider one example, if $P$ is
$\mathcal{B}(n)$ and $Q$ if the chain of length $n$ then the mapping from $P$ to $Q$ which sends each
set to its carcinality is order preserving. In this section we will see how an order preserving mapping
can be used to establish a relation between $\mu_P$ and $mu_Q$.

To begin, we introduce the M\"obius number of a poset. If $P$ is a poset, let $\widehat{P}$ be the poset obtained from $P$ by adjoining a new zero-element $\widehat{0}$ and a new one-element $\widehat{1}$. Hence if $x\in P$ then 
\[
	\widehat{0} \le_{\widehat{P}} x \le _{\widehat{P}} \widehat{1}.
\]
We define the M\"obius number $\mu(P)$ of $P$ by
\[
	\mu(P):=\mu_{\widehat{P}}(\widehat{0}, \widehat{1}).
\]
It is equal to the number of fchains of even length in $P$, less the number of chains of odd length. Note that the M\"obius number of the empty poset is $-1$ (Why?)

The following simple result will be one of our main tools. It implies that if a poset $P$
has a 1-element, then $\mu(P)=0$.

\begin{lemma}
	\label{lem:compmu0}
	If the poset $P$ has an element which is comparable with all elements of $P$ then $\mu(P)=0$.
\end{lemma}

\proof 
Suppose that $a$ is comparable with all eleents of $P$. Then there is bijection between the 
chains in $P$ which contain $a$ and whose which do not.\qed

If $a\in P$ and $a$ is comparable with every element of $P$, we will say that $P$ is 
a \textsl{cone} over $a$.

More definitions. Suppose that $P$ is a poset and $a\in P$. By $P_{a\le}$ we denote the set of 
elements $x$ of $p$ such that $a\le x$, while $P_{a<}$ consists of the elements $x$ such 
that $x>a$. Similarly we define $P_{\le a}$ and $P_{<a}$. Now we can state the main result of 
this section. It is more or less equivalent to Theorem 5.5 in Baclawski \cite{Baclawski1977}.

\begin{theorem}[Baclawski]
\label{thm:Bacl}
	Let $P$ and $Q$ be posets and let $f$ be an order-preserving mapping from $P$ to $Q$. Then 
	\[
		\mu(Q)=\mu(P)+\sum_{y\in Q} \mu(Q_{y<}) \mu(f^{-1}(Q_{\le y})).\qed 
	\]
\end{theorem}

A poset of the form $f^{-1}(Q_{\le t})$ will be called a \textsl{fibre} of $f$. Note that is a 
subset of $P$. The poset $Q_{\le y}$ has a $1$-element and so its M\"obius number is zero 
(by Lemma~\ref{lem:compmu0}). If all fibres of $f$ have $1$-elements then it follows from 
Theorem~\ref{thm:Bacl} that $\mu(P)=\mu(Q)$. A subset $P$ of a poset $s$ is an \textsl{ideal} if, whenever $a\in S$ and $x\le a$, we have $x\in P$. Any fibre of an order-preserving mapping is an ideal. It is worth noting that if $f$ is an order-preserving mapping from $S$ to the chain $\Ch(1)$ then $f^{-1}(0)$ is an ideal and, conversely, each ideal of $S$ determines an order-preserving mapping into $\Ch(1)$. A subset $F$ of $S$ is a \textsl{filter} if whenever $x\ge a$ and $a\in F$, we have $x\in F$.

The following result is a consequence of Theorem~\ref{thm:Bacl}, but we give a direct proof of it.

\begin{lemma}
\label{lem:Pideal}
	If $P$ is an ideal of the poset $S$ then 
	\[\mu(S) = \mu(P) + \sum_{y\in S\backslash P} \mu(S_{y<}) \mu(P_{\le y}).\]
\end{lemma}

\proof
We use Lemma~\ref{lem:Hall}. Suppose that $C$ is a chain in $S$. If $C\sbs P$ then, in the right side of map $a$, it is counted by the term $\mu(P)$. If $C\not\sbs P$, let $y$ be the least element of $C\backslash P$. Then $C\backslash y$ is the disjoint union of a chain from $P_{\le y}$ and a chain from $S_{y<}$. It is easy to check that in this case $C$ is counted, with the correct sign, by the expression $\mu(S_{y<}) \mu(P_{\le y})$.
\qed

Now we show how to derive Theorem~\ref{thm:Bacl} from Lemma~\ref{lem:Pideal}. Assume that $f$ is an order-preserving map from $P$ to $Q$. Construct a new poset $S$ with element set $P\cap Q$ by declaring that $a\le b$ if either
\begin{enumerate}[(a)]
	\item $a,b\in P$ and $a\le_P b$, or
	\item $a,b\in Q$ and $a\le_Q b$, or
	\item $a\in P, b\in Q$ and $f(a)\le_Q b$.
\end{enumerate}
(This construction is due to Baclawski \cite{Baclawski1977}.) It is easy to see that $S$ is a poset and $P$ is an ideal in it. Hence we have
\[\mu(S) = \mu(P)+\sum_{y\in S\backslash P} \mu(S_{y<})\mu(P_{\le y}).\]
If $y\in S\backslash P=Q$ then $S_{y<}=Q_{y\le }$ and $P_{y\le}=f^{-1}(Q_{\le y})$, whence we deduce that
\[\mu(S)=\mu(P)+\sum_{y\in Q} \mu(Q_{y<})\mu(f^{-1}(Q_{\le y})). \tag*{\sqr53}\]
(This is just a dual version of Lemma~\ref{lem:Pideal}.)

Since $P$ is an ideal in $S$ it follows that $Q^{\circ p}$ is an ideal in $S^{\circ p}$. Therefore
\begin{align*}
\mu(S) &= \mu(S^{\circ p})\\
&= \mu(Q^{\circ p}) + \sum_{x\in S^{\circ p}_{x<}}\mu(Q^{\circ p}_{\le x})\\
&= \mu(Q) + \sum_{x\in P} \mu(S_{<x})\mu(Q_{x\le}).
\end{align*}
As $Q_{x\le}=Q_{f(x)\le}$ has a $0$-element for all $x$ in $P$, this implies that $\mu(S)=\mu(Q)$. Hence Theorem~\ref{thm:Bacl} follows.

\textbf{Exercise}: Derive Lemma~\ref{lem:Pideal} from Theorem~\ref{thm:Bacl}.

In \cite{Walker1981} Walker proves a more general result than Theorem~\ref{thm:Bacl}: he allows the order-preserving mapping $f$ to be an \textsl{ideal relation} between $P$ and $Q$, i.e., an ideal in $P\times Q$. This has the advantages of being more general, and more symmetric in the roles $p$ and $Q$ play.

\section{Retracts}

A mapping $f:P\to P$ is \textsl{decreasing} if $f(x)\le x$ for all $x$ in $P$. A subset $Q$ of $S$ is a \textsl{retract} if there is an order-preserving and decreasing mapping $f$ from $S$ to $Q$ such that $f\upharpoonright Q$ is the identity, and then we call $f$ a \textsl{retraction}. If it is order-preserving and increasing then we also define the fixed points of $f$ to be a retract of $S$. Note that if $S$ is constructed from $P$ and $Q$ as in the proof of Theorem~\ref{thm:Bacl} then the mapping which sends $a$ in $Q$ to itself and $a$ in $P$ to $f(a)$ is a retraction. We saw in the proof of Theorem~\ref{thm:Bacl} that $\mu(S)=\mu(Q)$. More generally we have the following result.

\begin{lemma}
	If $Q$ is a retract of $S$ then $\mu(Q)=\mu(S)$.
\end{lemma}

\proof
Let $f$ be a retraction from $S$ onto $Q$. If $x\in f^{-1}(Q{\le y})$ then $f(x)\le y$. As $f(y)\le y$ (indeed $f(y) = y$) it follows that $y$ is a $1$-element $f^{-1}(Q_{\le y})$, hence each fibre of $f$ has a $1$-element and therefore has M\"obius number zero. By Theorem~\ref{thm:Bacl} we deduce that $\mu(Q) = \mu(S)$.\qed

\begin{corollary}
	If $P$ is a poset then $\mu(P) = \mu(\Ch(P))$.
\end{corollary}

\proof
Each element of $P$ is a chain, therefore $P$ is a subposet of $\Ch(P)$. Consider the map $f$ from $\Ch(P)$ to $P$ defined by setting $f(C)$ equal to $\max(C)$. Then $f$ is order-preserving, decreasing and its restriction to $P$ is the identity.
\qed

A \textsl{point} in a lattice is an element which covers $0$.

\begin{corollary}
	Let $P$ be a poset. If $1$ is not a join of points then $\mu_L(0,1)=0$.
\end{corollary}
\proof
If $x\in L\backslash 0$, define $f(x)$ to be join of the points in $L$ below $x$. Then $f$ is order-preserving and decreasing and $f(f(x))=f(x)$ for all $x \in L\backslash 0$. If $f(1)\ne 1$ then $f(1)<1$ and $F$ is retract of $L'$. Hence $\mu(F)=\mu(L')$ and, since $f(1)$ is a $1$-element in $F$, it follows that $\mu(L')=0$.
\qed

If $a$ and $b$ are elements of a poset $P$ and the least upper bound of $a$ and $n$ is defined, we denote it by $a\vee b$.

\begin{lemma}
\label{lem:avmu0}
	Let $P$ be a poset. If $a\in P$ and $a\vee x$ exists for all $x$ in $P$ then $\mu(P)=0$.
\end{lemma}

\proof
There are two steps. First, $a$ is a $0$-element in $P_{a\le}$ and so $\mu(P_{a\le})=0$. Second, the map $x\mapsto x\vee a$ is order preserving and increasing, with $P_{a \le}$ as its set of fixed points. Hence $P_{a \le }$ is a retract of $P$ and therefore $\mu(P)=0$.\qed

\begin{lemma}[Weisner]
	If $L$ is a lattice and $a\in L\backslash 0$ then
	\[
		\mu_L(0,1) = -\sum_{x\vee a=1, x<1} \mu_L(0,x).
	\]
\end{lemma}	

\proof
Suppose 
\[G := \{ x\in L': x\vee a<1\}.\]
Then $a\in G$ and $a\vee x$ exists for all $x$ in $G$, so $\mu(G)=0$ by the previous lemma. Since $G$ is an ideal in $L$, using Lemma~\ref{lem:muPab} we find that
\[\mu(G)=\sum_{x\in G\cup 0} \mu_L(0,x).\]
We now have
\begin{align*}
\mu_L(0,1)
&= -\sum_{x<1} \mu_L(0,x)\\
&= -\sum_{x\vee a=1, x<1} \mu_L(0,x) - \sum_{x\in G\cup 0} \mu_L(0,x)\\
&= -\sum_{x\vee a=1, x<1} \mu_L(0,x)-\mu(G)
\end{align*}
This yields the lemma.\qed

In delete, we will need the next result. The proof is left as an easy exercise.

\begin{lemma}
	Let $g$ be an order preserving and decreasing mapping of $P$ into itself and let $F$ be 
	the set of fixed points of $f$. Then $F$ is a retract of $P$.\qed
\end{lemma}

\section{Cutsets}

A \textsl{cutset} in a lattice $L$ is a set $C$ which contains at least one element from each 
maximal chain. Call a non-empty subset $S$ of $C$ a \textsl{simplex} if it has a bound 
(i.e., a meet or a join) in $L\backslash\{0,1\}$. The set of all simplices in $C$ forms a 
simplicial complex, which we denote by $\cS(L,C)$. By way of example, if $L$ is the lattice of 
subspace of a finite-dimensional vector space $V$ then the set of all $1$-dimensional vector spaces 
is a cutset, $C$ say. A subset of $C$ lies in $\cS(L,C)$ if and only if its join is not the entire 
space, i.e., if and only if it is not a spanning st in $V$. For any lattice $L$, let $L'$ denote 
the poset obtained from $L$ by deleting its $0$- and $1$-element. Thus
\[
	\mu(L') = \mu_L(0,1),
\]
which provides one reason why we need $L'$.

\begin{theorem}
	\label{thm:cuts}
	If $L$ is a lattice and $C$ is a cuset then $\mu(L')=\mu(\cS(L,C))$.
\end{theorem}

\proof
Let $\cS(L,C)$ be abbreviated to $\cS$. If $B$ is a chain in $L'$, define $f(B)$ by
\[
	f(B):=\{x\in C: x\cup B\in \Ch(L').
\]
In other words, $f(B)$ is the set of all elements of $C$ which are comparable with each element 
of $B$.) Since $C$ is a cutset, $f(B)\ne \emptyset$ and it follows that $f$ is an order-preserving mapping from $\Ch(L')$ to $?^{\op}$. Hence we may prove the theorem by showing that all fibres 
of $f$ have M\"obius number zero.

Let $S$ be an element of $\cS^{\op}$, and let $F$ denote the fibre $f^{-1}(\cS^{\op}_{\le S})$. If a chain $D$ of $L$ lies in this fibre then $s\sbs f(D)$. If $x\in D$ then $x$ is comparable with every element of $S$. Hence
\[
	\wedge S\le x \le \vee S.
\]
Since $S\in\cS(L,C)$, either $\wedge S$ or $\vee S$ lies in $L'$. Assume $\wedge S\in L'$, and denote it by $z$. Then for any element $x$ of $D$ we have $x\le z$, whence $D\cup z\in \Ch(L')$. Further every element of $S$ is comparable with all elements of $D\cup z$, thus $D\cup z$ belongs to $F$.

Now $z$ is a chain in $L'$ and $S\sbs f(z)$. Hence $z\in F$ and $z\cup D$ lies in $F$ for all elements of $F$. By Lemma~\ref{lem:avmu0} it follows that $F$ has M\"obius number zero. A similar argument yields the same conclusion when if $\wedge S\in L'$. Hence the theorem holds.\qed

We can manipulate this theorem to obtain a more explicit formula for $\mu(L')$.

\begin{corollary}
\label{cor:mucuts}
	Let $C$ be a cutset in the lattice $L$, and let $a_k$ be the number of $k$-subsets of $C$ 
	with join $1$ and meet $0$. Then 
	\[
		\mu_L(0,1)=\sum_{k=0}^C (-1)^ka_k.
	\]
\end{corollary}

\proof
If $S\in\cS$ then the interval $\hat{\cS}_{\le S}$ is a Boolean lattice and so
\[\mu_{\hat{\Gamma}}(0,S)=(-1)^{\abs{S}},\]
from we find, using Lemma~\ref{lem:muPab}, that
\begin{equation}
	\mu(\cS) = \mu_{\hat{\Gamma}(0,1)} 
		= \sum_{S\in \hat{\Gamma}\backslash 1} \mu_{\hat{\Gamma}}(0,S) 
		= 1+\sum_{S\in\Gamma}(-1)^{\abs{S}}.
\label{eqn9.1}
\end{equation}
Define $a_k$ to be the number of $k$-subsets $S$ of $C$ such that $\wedge S=0$ and $\vee S=1$. Then, when $1\le k\le \abs{C}$, the number of $k$-subsets of $?$ is equal to 
\[
	\binom{\abs{C}}{k}-a_k.
\]
Now, assuming $\abs{C}>0$ (which is the only interesting case)
\[
	0 = \sum_{S\sbs C, S\in \Gamma} (-1)^{\abs{S}} + \sum_S\sbs C, S\not\in \Gamma, (-1)^{\abs{S}}.
\]
From (\ref{eqn9.1}) we see that the first sum here is equal to $\mu(?)-1$, while the second sum 
equals $1+\sum_{k\ge 0}(-1)^k a_k$. The result follows.\qed

We consider applications of Corollary~\ref{cor:mucuts}. Suppose that $L$ is the lattice of subspace of a finite-dimensional vector space over some finite field, and let $C$ be the set of all $1$-dimensional subspaces. Then $C$ is a cutset and $a_k$ is the number of spanning subsets of $C$ with cardinality $k$.

For another example, let $G$ be a graph with vertex set $V$ and let $L$ be the set of all partitions of $V$ such that each cell induces a connected subgraph of $G$. Then the join of any two elements of $L$ lies in $L$ and hence $L$ is a lattice, but not in general a sub-lattice of a lattice of all partitions of $V$, Let $C$ be the set of all partitions in $L$ with one cell of size two and all others singletons. (So the cell of size two is an edge of $G$.) Then $C$ is a cutset in $L$ and $a_k$ is the number of subgraphs of $G$ with $k$ edges and the same number of connected components as $G$.

Rota \cite{Rota1964} proved Theorem~\ref{thm:cuts} under the assumption that $C$ was a cutset and an antichain. Walker proves an even more general result than Theorem~\ref{thm:cuts} in \cite{Walker1981}, out proof is based on his. (Our task is slightly more complicated, in that Walker can use ideal relations where we must use order-preserving functions.)

\section{Complements}
The main result in this section is a slight weakening of Theorem~8.1 from Walker \cite{Walker1981}. 
To begin, we derive a technical lemma.

\begin{lemma}
\label{lem:mucomp}
	Let $L$ be a lattice and suppose $s\in L'$. Let $G$ be the set of all elements $x$ of $L'$ 
	such that $x\vee s<1$. If $\mu(G_{\le y})\ne 0$ then $y$ is a complement to $s$.\qed
\end{lemma}

\proof
Note that $s\in G$. We show that if $y$ is not a complement to $s$ in $L$ then $\mu(G_{\le y}) =0$. 
If $y\in G$ then $G_{\le y}$ has a $1$-element and so its M\"Obius number is zero. If $y\notin G$ 
and $y\wedge s=0$ then $y$ is a complement to $S$. Suppose $y\notin G$ and $t\wedge s \ne 0$. 
If $z\in G_{\le y}$, then
\[
	z\vee (y\wedge s) \vee s= z\vee s <1.
\]
Hence $z\vee (y\wedge s)\in G$ and therefore $z\in G_{\le y}$. Thus $G_{\le y}$ is a cone over 
$y\wedge s$ and so has M\"obius number zero. The lemma follows.\qed

\begin{theorem}[Walker \cite{Walker1981}] 
\label{thm:Walker}
	Let $L$ be a lattice, let $a$ be an element of $L'$ and let $a^-$ be the set of all complements 
	of $a$. Then $\mu(L'\backslash a^-)=0$.
\end{theorem}

\proof
Let $M$ denote $L'\backslash a^-$ and let $G$ be the subposet of $L'$ consisting of all elements $x$ of $M$ such that $x\vee a<1$. Then $G$ is an ideal of $L'$ and contains $a$. The fibres of the inclusion mapping of $G$ in $M$ are the sets $G_{\le y}$, where $y$ in $M$. By the previous lemma these fibres all have M\"obius number zero, whence $\mu(M)=\mu(G)$. Since $a\in G$ and $a\vee x$ exists for all $x$ in $G$, we see by Lemma~\ref{lem:avmu0} that $\mu(G)=0$. It follows that $\mu(M)=0$, as required.\qed

\textbf{Exercise}: Let $s$ be an element of the lattices $L$. Let $S$ be a subset of $L'$ containing $a^-$ such that if $x\in S$ then $x\vee a=1$. Show that $\mu(L'\backslash S)=0$.

\begin{corollary}
	If $L$ is a lattice and $\mu_L(0,1)\ne 0$ then $L$ is complemented.
\end{corollary}

\proof
If some element of $L$ has no complement then the theorem applies, with $S\ne \emptyset$.\qed

\section{Topology}

There is more going on than we have yet admitted. An order-preserving map $f$ from a poset $P$ to a poset $Q$ induces an order-preserving map from $\Ch(P)$ to $\Ch(Q)$. But $\Ch(P)$ and $\Ch(Q)$ are simplicial complexes and thus may be viewed as topological spaces. The map induced by $f$ is then a continuous map.

Let $X$ and $Y$ be topological spaces. Two continuous functions $f$ and $g$ from $X$ to $Y$ are \textsl{homotopic} if there is a continuous function
\[
	\Phi: X \times [0,1] \to Y
\]
such that $\Phi(x,0)=f(x)$ and $\Phi(x,1)$ for all $x$ in $X$. We say two topological spaces $X$ 
and $Y$ are homotopic if there are continujous functions $f:X\to Y$ and $g:X\to Y$ such 
that $g\circ f$ and $f\circ g$ are homotopic to the respective identity maps on $X$ and $Y$. 
It can be shown that homotopy is an equivalence relation on topological spaces. Any convex subset 
of $\re^n$ is homotopic to a point, while two homotopic surfaces in $\re^3$ are homeomorphic. 
We say posets $P$ and $Q$ are homotopic if $\Ch(P)$ and $\Ch(Q)$ are.

For our purposes, the following is important.

\begin{lemma}
	If $P$ and $Q$ are posets such that $\Ch(P)$ is homotopic to $\Ch(Q)$ then $\mu(P)=\mu(Q)$.
\end{lemma}

\proof
The M\"obius number of $P$ is determined by the Euler characteristic of $\Ch(P)$. Homotopic simplicial complexes have the same Euler characteristic.\qed

A topological space is \textsl{contractible} if it is homotopic to a point. One class of contractible simplicial complexes are \textsl{cones}. A simplicial complex $\mathcal{S}$ is a cone if it contains an element $v$ such that $v\vee x$ is defined for all elements $x$ of $\mathcal{S}$. It is not hard to see that if the poset $P$ is a cone then $\Ch(P)$ is a cone as a simplicial complex (and as a poset). We will say a poset $P$ is contractible if $\Ch(P)$ is. We have the following important result.

\begin{theorem}[Quillen \cite{Quillen1978}]
	Let $P$ and $Q$ be posets and let $f$ be an order-preserving map from $P$ to $Q$. If $f^{-1}(q)$ 
	is	contractible for any element $q$ of $Q$ then $P$ and $Q$ are homotopic.\qed
\end{theorem}

Note that in all cases where we have proved that the fibres of some order-preserving map 
have M\"obius number zero, we have actually shown that the fibres are cones and hence contractible. 
Thus if $L$ is a lattice, $p\in L$ and $S$ is the set of complements of $s$ in $L$, 
then $L'\backslash S$ is contractible.

If $P$ and $Q$ are posets, it makes sense to talk about two order-preserving maps $f$ and $g$ from $P$ to $Q$ as being homotopic. No combinatorial characterization of what this means is known. However if $f(x)\le g(x)$ for all $x$ in $P$ then it is easy to show that $f$ is homotopic to $g$.

We may view $Q^P$ as a poset, where $f\le g$ for two elements $f$ and $g$ of $Q^P$ if $f(x)\le g(x)$ for all $x$ in $P$. The poset $Q^P$ is the disjoint union of a number of connected components; two maps in the same component will be homotopic. The constant map taking each element of $P$ to a fixed element of $Q$ is always order-preserving, so $\abs{Q^P}\ge \abs{Q}$.

It turns out that lying in the same component of $Q^P$ is not a good approximation to the topological notion of homotopy, for reasons we now discuss.

Suppose that $a$ is an element of $P$ which covers a unique element $b$ of $P$. Define a map 
\[
	\phi_q : P\to P\backslash a
\]
by setting $\phi_a(x)$ equal to $x$ if $x\ne a$ and $\phi_a(a)=b$. Call $\phi_a$ a \textsl{deletion}. Each fibre of the inclusion mapping of $P\backslash a$ into $P$ has a $1$-element. (It is not hard to see that $P\backslash a$ is a retract of $P$.)

Now suppose that $f$ is an order-preserving map of $P$ into itself and $f(x)\le x$ for all elements $x$ of $P$. Let $a$ be an element of $P$ that is minimal, subject to the condition that $f(a)<a$. If $b<a$ then, by our choice of $a$, we have $f(b)=b$. On the other hand, $f$ is order-preserving and so $f(b)\le f(a)$. Hence if $b<a$ then $f(b)\le f(a)$ and we have shown that $f(a)$ is the unique element of $p$ covered by $a$.

\textbf{Exercise}: Show that any order-preserving and decreasing map from $P$ into itself is a composition of deletions.

It is possible that, by applying a sequence of deletions, we might be able to map $P$ onto the poset with exactly one element. In this case we say that $P$ is \textsl{dismantlable}, and $P$ is homotopic to a $1$-element poset, i.e., it is contractible.

\textbf{Exercise}: Let $P$ be a poset. The following are equivalent:
\begin{enumerate}[(a)]
\item $P$ is dismantlable,
\item $P^X$ has exactly one component for any poset $X$ a nd
\item $P^P$ has exactly one component.
\end{enumerate}

\section{Geometric Lattices}

The M\"obius function is particularly useful when applied to geometric lattices. This section introduces these lattices briefly. There are two parts to their definition.

A lattice is a \textsl{point-lattice} if every non-zero element can be expressed as the join of points.

A lattice is \textsl{semimodular} if, whenever $a$ and $b$ are elements such that if $a$ covers $a\wedge b$ then $b$ is covered by $a\vee b$. There are a number of equivalent definitions, which we will discuss shortly. What we have just callsed a semimodular lattice is more strictly an upper semimodular lattice. A lattice which is dual to a semimodular lattice is \textsl{lower semimodular}.

A lattice is \textsl{geometric} if it is a semimodular point lattice. One class of examples arises 
as follows. Let $X$ be a set of points in a finite-dimensional projective space and define 
a \textsl{flat} in $X$ to bw any subset of $X$ of the form $H\cap X$, where $H$ is a 
projective subspace. The lattice of flats of $X$ is geometric.

\textbf{Exercise}: Show that $\mathcal{P}(n)$, the lattice of all partitions of an $n$-set, is geometric.

The points of a geometric lattice may also be referred to as \textsl{atoms}. A maximal flat is called a \textsl{hyperplane}.

\textbf{Exercise}: Show that each element in a geometric lattice is the meet of a set of hyperplanes.

We will use the result of the next exercise several times.

\textbf{Exercise}: Show that any interval in a geometric lattice is geometric.

For the remainder of this section, we discuss some of the properties of semimodular lattices.

\begin{lemma}
A lattice $L$ is semimodular if $a\vee b$ covers both $a$ and $b$ whenever $a$ and $b$ cover $a\wedge b$.
\qed
\end{lemma}

A poset $P$ is \textsl{ranked} if any two maximal chains joining the same pair of elements have the same length. (Equivalently, we may say that $P$ satisfies the Jordan-Dedekind condition.) If $P$ is ranked and $a\in P$ then the maximum length of a chain ending on $a$ is the rank of $a$, which we denote by $r(a)$. If $P$ is a ranked poset and $b$ covers $a$ then $r(b)=r(a)+1$.

\begin{lemma}
Let $L$ be a lattice with rank function $f$. Then $L$ is semimodular if and only if 
\begin{equation}
r(a\wedge b) + r(a\vee b) \le r(a) + r(b)
\label{12.1}
\end{equation}
\end{lemma}

We call Equation~(\ref{12.1}) the \textsl{semimodular identity}. A ranked lattice is modular if equality holds in the semimodular identity for all pairs of elements $a$ and $b$. The lattice of subspaces of a vector space is modular, as are the Boolean lattices.

We will need the next result in Section~\ref{sec:geomoeb}.

\begin{lemma}\label{12.3}
Any point in a geometric lattice has a complement.
\end{lemma}

\proof
Let $L$ be geometric and lert $p$ be a point in $L$. Let $a$ be an element of $L$ which is maximal, subject to the condition that $a\wedge p=0$. If $a\vee p=1$ then $a$ is a complement of $p$ and we are finished.

Otherwise $a\vee p<1$ and, since $1$ is join of points, it follows that there is point $q$ 
of $L$ such that $q\not \le a\vee p$. Now there are two possibilities. If $p\le a\vee q$,
then 
\[
	p\vee a\le q\vee a.
\]
But $p$ and $q$ cover $0$, hence both $a\vee p$ and $a\vee q$ cover $a$. This implies 
that $r(a\vee p)=r(a\vee q)$ and therefore $a\le a\vee p$, which contradicts our choice of $q$.

If $p \not \le a\vee q$ then $p\wedge (a\vee q)=0$. Since $q\not \le a\vee p$, it follows 
that $a<a\vee q$, and thus we have a contradiction to our choice of $a$.\qed

For further background on geometric lattices, see \cite{Baclawski1977,Rota1964}.

\section{Modular Elements}

An element $a$ in a geometric lattice $L$ is \textsl{modular} if the semimodular identity holds 
for all pairs $(a,b)$, i.e.,
\[
	r(a\wedge b) + r(a\vee b) = r(a)+ r(b)
\]
for all elements $b$ of $L$. Equivalently $a$ is modular if the set of all complements of $L$ 
form an antichain.

\textbf{Exercise}: Prove that an element $a$ in a geometric lattice is modular if and only if its complements form an antichain.

Any point in a geometric lattice $L$ is modular. If $a$ is a point of $L$ and $b\in L$ then, since $a$ covers $0$, either $a\wedge b=0$ or $a\wedge b=a$. In the firs case
\[r(x) +1 = r(x) + r(a) \ge r(x\wedge a) + r(x\vee a) = r(x\vee a) > r(x),\]
while in the second
\[r(x) + r(a) \ge r(x\wedge a) + r(x\vee a) =r(a) + r(x).\]
In both cases we have equality in the semimodular inequality. The fact that points are modular is not always useful; it may be better to have modular elements of higher rank. We note two examples.

Let $\mathcal{B}_q(n)$ denote the lattice of subspaces of an $n$-dimensional vector space over a field with $q$ elements. This is a modular lattice, and thus all tis elements are modular. If $\mathcal{P}(n)$ is the lattice of all partitions of $\{1,\cdots,n\}$ then the partition with cells $\{1,\cdots,n-1\}$ and $\{n\}$ is a modular hyperplane.

\begin{lemma}\label{13.1}
Let $a$ and $b$ be elements in the geometric lattice $L$. If $a$ is modular then the 
map $x\mapsto x\vee a$ is an isomorphism from $[a\wedge b, a]$ to $[b,a\vee b]$.
\end{lemma}

\proof
If $x\in[a\wedge b, a]$ then the mapping $x\mapsto b\wedge x$ is order-preserving, as is the 
mapping $y\mapsto a\vee y$ when $y\in [a, a\vee b]$. Hence the composite map $\psi$ defined by
\[
	\psi(x) = a\wedge (x\vee b)
\]
is an order-preserving map from $[a\wedge b, a]$ into itself. Since $a\wedge(x\vee b) \ge x$, 
it is also increasing.

Suppose $c\in [a\wedge b,a]$. Since $a\wedge b = c\wedge b$, the semimodular identity implies that
\[
	r(c\vee b) - r(b) = r(c) - r(a\wedge b).
\]
Applying the semimodular identity to the pair $(a, \vee b)$ and noting that 
$(a\vee c) \vee b = b\vee a$, we get
\[
	r(a\vee b) - r(c\vee b) \le r(a) - r(a\wedge (c\vee b)).
\]
Summing the last two inequalities yields that
\[
	r(a\vee b) - r(b) \le r(a) - r(a\wedge b) - (r(a\wedge(c\vee v)) - r(c)).
\]
As $a$ is modular
\[
	r(a\vee b) - r(b) = r(a) -r(a\wedge b)
\]
and, given the previous inequality, we deduce that
\[
	r(a\wedge (c\vee b)) \le r(c).
\]
However $c\le a\wedge (c\vee b)$ and therefore we have proved that $c=a\wedge (c\vee b)$. So 
if $a$ is modular then $\psi$ is the identity mapping andthe intervals $[a\wedge b, a]$ 
and $[v, a\vee b]$ are isomorphic.\qed

\section{M\"obius Functions and Geometric Lattices}
\label{sec:geomoeb}

Our first result will enable us to compute the M\"obius function on intervals in $\mathcal{B}_q(n)$ 
and $\mathcal{P}(n)$. We need one preliminary result.

\begin{lemma}\label{14.1}
If $C$ is an antichain in the poset $P$ then
\[
	\mu(P) = \mu(P\backslash C) + \sum_{x\in C} \mu(P_{a_<}\mu(P_{<a}).
\]
\end{lemma}

\proof
Apply Theorem~\ref{thm:Bacl} with $f$ the inclusion mapping of $P\backslash C$ into $P$. The details 
are left as an exercise.\qed

This lemma is useful even when $C$ is a single element of $P$.

\begin{theorem}\label{14.2}
Suppose $a$ is a modular element of the geometric lattice $L$, not $0$ or $1$, and let $a^-$ be 
the set of all complements of $a$ in $L$. Then
\[
	\mu_L(0,1) = \mu_L(0,a) \sum_{x\in a^{\perp}} \mu_L(0, x).
\]
\end{theorem}

\proof
Let $P$ be $L'$ and let $a^-$ be the set of all complements to $a$ in $L$. Then $a^-$ is an antichain and so, using Lemma~\ref{14.1}, we get
\[
	\mu(L') = \mu(L'\backslash a^-) + \sum_{x\in a^{\perp}} \mu(L'_{x<}) \mu(L'_{<x}).
\]
By Theorem \ref{thm:Walker} we have that $\mu(L'\backslash a^-)=0$. By Lemma \ref{13.1}, if $x\in a^-$ then the intervals $[x\wedge a,a]=[0,a]$ and $[x, x\vee a]=[x,1]$ are isomorphic, hence
\[
	\mu(L'_{x<}) = \mu_L(x,1) = \mu_L(0,a).
\]
As $\mu(L'_{<x}) = \mu_L(0,x)$, the theorem follows.\qed

\begin{corollary}
	If $L$ is a geometric lattice and $a$ and $b$ are elements of $L$ such that $a\le b$ then 	$(-1)^{r(b)-r(a)} \mu_L(a,b)>0$.
\end{corollary}

\proof
Since any interval of a geometric lattice is geometric, it suffices to assume that $a=0$ and $b=1$. Let $p$ be a point in $L$. Then $p$ is modular and all its complements are hyperplanes. (It has complements by Lemma~\ref{12.3}.) By the theorem
\[
	\mu_L(0,1) = -\sum_{x\in p^{\perp}} \mu(0,x).
\]
We may assume inductively that $\mu(0,x)$ is non-zero and has the same sign for all $x$ in $p^-$, whence the result follows.\qed

Next we compute the M\"obius function on $\mathcal{B}_q(n)$ using Theorem \ref{14.2}. Let $h$ 
be a hyperplane in $\mathcal{B}_q(n)$. Then $h$ is modular and so, if $L=\mathcal{B}_q(n)$,
\[
	\mu_L(0,a) = \mu_L(0,h) \sum_{p\in h^{\perp}} \mu_L(0,p).
\]
Since $h$ is modular, all its complements are points. Consequently $\mu_L(0,p)=-1$. The number of points in $h^-$ is $q^{n-1}$ and therefore
\[
	\mu_L(0,1) = -q^{n-1} \mu_L(0,h).
\]
As $\mu_L(0,h) = \mu_{\mathcal{B}_q(n-1)} (0,1)$, a trivial induction argument yields that
\[
	\mu_{\mathcal{B}_q(n)}(0,1) = (-1)^n q^{\binom{n}{2}}.
\]

We can also compute the M\"obius function for $\mathcal{P}(n)$. Here $h=\{\{1\}, \{2,\cdots.n\}\}$ is a modular hyperplane whose complements are the partitions with one non-trivial cell, of the form $\{1,i\}$. Hence
\[
	\mu_{\mathcal{P}(n)} (0,1) = -(n-1) \mu_{\mathcal{P}(n-1)} (0,1).
\]
Once again a simple induction argument yields that
\[
	\mu_{\mathcal{P}(n)}(0,1) = (-1)^{n-1} (n-1)!.
\]

\section{Broken Circuits}
The main result of this section shows that if $L$ is a geometric lattice then $(-1)^r\mu_L(0,1)$ 
is not only non-negative, it counts something.

Let $L$ be a geometric lattice and let $S$ be the set of all points in it. Since $L$ is 
a point-lattice, we can identify each element of $L$ with the set of points below it in $L$. 
We can extend the rank function of $L$ to a function on subsets of $S$ by defining $r(T)$ 
to be $r(\vee T)$, for any subset $T$ of $S$. We have

\begin{enumerate}[(1)]
	\item $r(\emptyset)=0$,
	\item if $p\in S$ then $r(p)=1$,
	\item if $T$ and $U$ are subsets of $S$ and $T\sbs U$ then $r(T)\le r(U)$ and 
	\item for any pair of subsets $T$ and $U$ of $S$,
		\[
			r(T) + r(U) \ge r(T\cup U) + r(T\cap U).
		\]
\end{enumerate}

We define a subset $T$ of $S$ to be \textsl{independent} if $r(T)=\abs{T}$, all other subsets are \textsl{dependent}. The set $S$, together with its collection of independent subsets, is a \textsl{matroid}. A \textsl{circuit} is a minimal dependent subset of $S$. A \textsl{flat} is a subset, $F$ say, of $S$ such that if $p\in S\backslash F$ then $r(\cup F)>r(F)$. Thus the flats correspond precisely to the elements of $L$. We will not be doing any matroid theory, but we  will need to refer to the circuits and independent sets of a geometric lattices.

The independent sets of a geometric lattice form a simplicial complex--every subset of an independent set is independent. We are now going to define the \textsl{broken circuit} complex, which is a subcomplex of the independent set complex. Assume $L$ is a geometric lattice and let $\unlhd$ be a total order on its points. A set of points is a \textsl{broken circuit} if it can be obtained from some circuit by deleting its least element, relative to $\unlhd$. The broken circuit complex $\Br(L)$ has as its elements all independent sets which do not contain a broken circuit. Since a set of points with contains no broken circuit cannot contain a circuit, then elements of $\Br(L)$ are all independent sets. If $T\in \Br(L)$ then $\vee T$ contains no point less than the least element of $T$. (If $p$ is a point in $\vee T$ and $p\notin T$ then there is a circuit $p\cup S$, for some subset $S$ of $T$.)

\begin{theorem}[Whitney \cite{Whitney1932}]\label{15.1}
	Let $L$ be a geometric lattice and let $\unlhd$ be a total order on its points. Then 
	the number of independent sets of $k$ points which contain no broken circuit 
	is $\sum_{a: r(a)=k} (-1)^k \mu_L(0,a)$.
\end{theorem}

\proof
Any independent set of size $k$ lies in a unique element of $L$ with height $k$. Hence it suffices to prove that $(-1)^r \mu_L(0,1)$ is the number of independent sets of $r$ atoms containing no broken circuits, where $r$ is the height of $L$. We prove this by induction on $r$.

Let $p_1$ be the least point of $L$ and let $bn$ be a complement of $p$ in $L$. Then $b$ has height $r-1$. Let $M$ the geometric lattice formed by the interval $[0,b]$ and let $\Br(M)$ be the broken circuit complex of $M$, relative to the ordering of the points of $M$ obtained by restriction of $\unlhd$. We claim that $T$ is an independent set of $r-1$ points of $M$ then $T\cup p_1 \in \Br(L)$ if and only if $T\in \Br(M)$.

Suppose first that $T\in \Br(M)$. If $T\cup p_1$ contains a broken circuit $C$ from $L$ then either $C\sbs T$ and so $C$ is a broken circuit in $M$, or $C\cup p_1$ is a circuit in $L$ and therefore 
\[
	p_1 \in \vee C\le b.
\]
Conversely, let $S$ be an $r$-subset in $\Br(L)$. Then all points of $L$ lie in $\vee S$, and therefore $p_1\cup S$ is dependent. Since $S$ is independent any circuit in $p_1\cup S$ must contain $p_1$, whence $S$ contains a broken circuit.

By induction, the number of $(r-1)$-subsets of $\Br(M)$ is equal to
\[
	(-1)^{r-1} \mu_M(0,1) = (-1)^{r-1} \mu_L(0,b)
\]
and therefore the number of $r$-sets in $\Br(L)$ is equal to
\[
	(-1)^{r-1} \sum_{b\in p_1^{\perp}} \mu_L(0,b).
\]
Since $\mu_L(0,p_1)=-1$, by~Theorem \ref{14.2} this last sum is equal to $(-1)^r \mu_L(0,1)$, 
as required.\qed

\textbf{Exercise}: A simplicial complex is \textsl{pure} if all its maximal elements have the same height. Show that any broken circuit complex is pure.

\section{The Partition Lattice}
In this section we apply Theorem \ref{15.1} to the partition lattice $\mathcal{P}(n)$. We identify the points of $\mathcal{P}(n)$ with the edge set of $K_n$ and we let $\unlhd$ denote the lexicographic order on the points. An independent set is then a forest, i.e., an acyclic subgraph of $K_n$. If $i<j<k$ then the edges $ik$ and $jk$ form a broken circuit. It follows that a forest $F$ in $E(K_n)$ contains no broken circuit if and only if each component of $F$ has the property that the vertices in any path going away from the least vertex form an increasing sequence. (This condition is equivalent to containing no ``broken triangle", the details are up to you.) Consequently the forests in $K_n$ containing no broken circuits can be viewed as non-increasing functions on the set $\{1,\cdots,n\}$, the number of components in the forest is equal to the number of fixed points of the function.

Since $1$ is a fixed point of any non-increasing function, the number of such functions with exactly one fixed point is $(n-1)!$. This shows that
\[
	\mu_{\mathcal{P}(n)}(0,1) = (-1)^{n-1} (n-1)!.
\]
Fortunately this is consistent with our earlier result. The problem which remains is to determine the number of non=increasing functions with exactly $k$ fixed points when $k>1$. I claim that this is equal to the number of permutations of $\{1,\cdots,n\}$ with exactly $k$ cycles.

The proof of this is indirect. The first step is an encoding of a permutation in cyclic written out explicitly. (To give an extreme case, the identity permutation in cyclic form is usually written as $(1)$, but we must write it as $(1)(2)\cdots(n)$.) Now write each cycle so that the largest element is first (so if $n=4$ then $(123)$ is now $(312)(4)$). Next, order the cycles so that the first elements form an increasing sequence. Finally remove the parentheses. You are invited to prove that we have now defined a bijection from $\sym{n}$ onto itself. Denote the image of a permutation $\beta$ under this bijection by $\hat{\beta}$.

What is the relation between the cycles of $\beta$ and $\hat{\beta}$? If $\sigma\in \sym(n)$, define $j$ to be a \textsl{record} if $\sigma(i)<\sigma(j)$ whenever $i<j$. Our claim (well it is my claim, but you have to prove it) is that the number of cycles in $\beta$ is equal to the number of records of $\hat{\beta}$.

But this only completes the first step; we need to convert $\hat{\beta}$ into a non-decreasing function. This is easy. If $\sigma\in \sym{n}$, let $f_{\sigma}$ be defined by 
\[
	f_{\sigma}(j) = \abs{\{i: i<\sigma^{-1}(j), \sigma(i) \le j\}}.
\]
Again, you must convince yourself that $\sigma$ can be reconstructed from $f_{\sigma}$. Note however that the fixed points of $f_{\sigma}$ are precisely the records of $\sigma$. Hence the number of fixed points of $f_{\hat{\beta}}$ is equal to the number of cycles of $\beta$.

A forest in $K_n$ with exactly $n-d$ edges has exactly $d$ components. So we have shown that the number of forests with $k$ edges which contain no broken circuit is equal to the number of permutations of $\{1,\cdots,n\}$ with exactly $n-k$ cycles. This number has no nice explicit form, but it is known to be equal to $(-1)^k$ times the coefficient of $x^k$ in $x(x-1)\cdots(x-n+1)$. The coefficient itself is a \textsl{Stirling number of the first kind}. (For background see, e.g., \cite[Ch.~1]{Stanley2002}.)

\section{Contractions and Colourings}
We consider a family of geometric lattices including $\mathcal{P}(n)$ as a special case. Let $G$ be a graph with vertex set $V$ and edges set $E$. A \textsl{contraction} of $G$ will be defined to be a partition of $V$ such that the subgraph induced by any cell is connected. Equivalently we may view them as subsets of $S$ of $E$ with the property that, for any edge $f\in E\backslash S$, the number of components of $S\cup f$ is less than the number of components of $S$. We will denote the set of all contractions of $G$ by $L_G$. Every contraction of $G$ is a partition of $V$ and so, if $n=\abs{V}$, it can be viewed as an element of $\mathcal{P}(n)$. Further the join of any two contractions is a contraction, and so the contractions of $G$ form a sub-semilattice of $\mathcal{P}(n)$. As we remark in the Appendix, any join-semilattice with zero can be turned into a lattice--in this case we define the meet of contractions $\sigma$ and $\tau$ by
\[
	\sigma \wedge \tau:= \vee \{\gamma \in L_G: \gamma \le \sigma, \tau\}.
\]

\textbf{Exercise}: Show that $L_G$, as defined above is a geometric lattice, and that if $S$ is a set of points of $L_G$ then $n-r(S)$ is the number of components of $S$.

The points of $L_G$ are precisely the edges of $G$. The independent sets of points are precisely the (edge-sets of the) forests in $G$, and the circuit are the circuits. The main result of this section will be an expression for the number of proper $k$-colourings of $G$ in terms of the M\"obius function of $L_G$.

A \textsl{proper $k$-colouring} of $G$ is a mapping $f: V\to \{1,\cdots,k\}$ such that $f(u)\ne f(v)$ whenever $uv$ is an edge. If $f$ is a mapping from $V$ to $\{1,\cdots,k\}$, define the set $K(f)$ by 
\[K(f):=\{uv\in E: f(u) = f(v)\}.\]
Note that the components of $K(f)$ form a contraction and that $f$ is a proper colouring of $G$ if and only if $K(f)=\emptyset$. Let $F_=(A,k)$ denote the number of mappings $f$ from $V$ to $\{1,\cdots,k\}$ such that $K(f)=A$ and let $F_{\le} (A,k)$ be the number of mappings from $V$ to $\{1,\cdots,k\}$ such that $K(f)\ge A$. Since $B\sbs K(f)$ if and only if $f$ is constant on the components of $B$, we have
\[F_{\le}(B,k) = k^{n-r(B)}.\]
As we also have
\[F_{\le}(B,k) = \sum_{A\sps B} F_= (A,k),\]
by M\"obius inversion we find that
\[F_=(A,k) = \sum_{B: B\sps A} \mu_{L_G}(A,B) k^{n-r(B)}.\]
Thus we have proved the following.

\begin{theorem}
	The number of proper $k$-colourings of the graph $G$ is equal to
	\[
		\sum_i \left(\sum_{\abs{A}=i} \mu_{L_G}(0,A)\right) k^{n-i}.\qed
	\]
\end{theorem}

In other words
\[
	\sum_i \left(\sum_{\abs{A}=i} \mu_{L_G}(0,A)\right)x^{n-i}
\]
is the chromatic polynomial of $G$. We have shown that the coefficients of the chromatic polynomial 
are the level numbers of the broken circuit complex of $G$. (Perhaps we should say of $L_G$.)

Let $\omega$ be a set of points in the projective space $PG(d,q)$. Then, as we noted earlier, the intersections of $S$ with the hyperplanes of $PG(n,q)$ are the elements of a geometric lattice. The rank of a subset of $\Omega$ is equal to the dimension of the space spanned by it. We are interested in counting the number of hyperplanes of $PG(n,q)$ that contain no point of $\Omega$. (Well, I am interested. You may have to fake it.)

The points of $\Omega$ can be represented by vectors $\seq{x}{1}{2}{n}$ in $V(d+1,q)$. If $a\in V(d+1,q)$ then the vectors $x_i$ such that $a^Tx_i=0$ are a hyperplane in the geometric lattice $L$ determined by $\Omega$. Denote the hyperplane corresponding to the vector $a$ by $h(a)$. If $S\sbs \Omega$, define $f(S)$ to be the number of vectors $a$ such that $h(a)=S$ and let $g(S)$ be the number of vectors $a$ such that $h(a)\sps S$. Then $g(S)=q^{d-r(S)}$ and consequently
\[
	f(S) = \sum_T \mu_L(S, T)q^{d-r(T)}.
\]
Therefore
\[
	f(\emptyset) = \sum_T \mu_L(0,T) q^{d-r(T)} 
		= \sum_{k\ge 0} \left(\sum_{r(T)=k} \mu_L(0,T)\right) q^{d-k}.\]
This is a polynomial in $q$, which we will denote by $F_L(q)$. It is called 
the \textsl{characteristic polynomial} of $L$.

\textbf{Exercise}: Show that the number of $t$-tuples of vectors $\seq{a}{1}{2}{t}$ such that $\cap_i h(a_i)=\emptyset$ is equal to $F_L(q^t)$.

There is coding theory view of all this, which is both interesting and useful. Suppose that we arrange the vectors $\seq{x}{1}{2}{n}$, given above into a $(d+1)\times n$ matrix, $G$ say. The row space of $G$ is a linear code over the field $GF(q)$. If $a\in V$ then $a^TG$ is a code word and the weight of this word is the number of elements of $\Omega$ not in $h(a)$. Thus the hyperplanes of the lattice of flats of $\Omega$ correspond to the code words with minimal non-zero weight. Further, there is a vector $a$ such that $h(a)$ is disjoint from $\Omega$ if and only if there is a code word with weight $n$, and the number of such codewords is equal to $f(\emptyset)$.

\section{Points and Hyperplanes}
\label{sec:py-hyps}

The main result in this section is that a geometric lattice always has at least as many hyperplanes as points. The lattice of subspaces of a finite vector space shows that equality can occur. The proof makes use of another interesting result.

\begin{theorem}[Dowling and Wilson \cite{DowlingWilson}]
\label{18.1}
	Let $L$ be a finite lattice. If $\mu_L(p,1)\ne 0$ for all elements $p$ of $L$ then there 
	is a permutation $\sigma$ of the elements of $L$ such that $q\vee \sigma(q)=1$ for all $q$ 
	in $L$.
\end{theorem}

\proof
We use Lemma~\ref{lem:gxysum}. Let $g$ be the real-valued function on $L$ defined by
\[g(p) = \begin{cases}
1,\quad \text{if } p=1;\\
0,\quad \text{otherwise}.
\end{cases}\]
Let $G$ be the matrix with rows and columns indexed by the elements of $L$ and with $(G)_{pq} = g(p\vee q)$. We can complete the proof by showing that $\det(G)\ne 0$. By Equation (1) from dets we have
\[\det(G) = \prod_x \sum_y \mu_L(x,y) g(y)=\prod_x \mu_L(x,1).\]
and, by our hypothesis on $L$, it follows that $\det H\ne 0$.\qed

Any permutation $\sigma$ satisfying the condition of Theorem \ref{18.1} must map $0$ to $1$. Hence if $L$ is geometric and $p$ is a point of $L$ then $\sigma(p)$ must be a hyperplane. Therefore $\sigma$ determines an injection of the points of $L$ into its hyperplanes, and so the number of hyperplanes in a geometric lattice is at least as large as the number of points. Actually a somewhat stronger statement can be made.

If $L$ is a lattice let $W_k$ denote the number of elem,emnts in $L$ with height $k$. If $L$ is geometric with height $n$ then $W_0=W_n=1$, while $W_1$ is the number of points and $W_{n-1}$ is the number of hyperplanes. We have just seen that Theorem \ref{18.1} implies that $W_1\le W_{n-1}$. (The numbers $W_i$ are sometimes referred to as the \textsl{Whitney numbers of the first kind}.)

\begin{corollary}\label{18.2}
	If $L$ is a geometric lattice with rank $d$ then
	\[
		W_0+\cdots + W_k \le W_{d-k} + \cdots + W_d.
	\]
\end{corollary}

\proof
Assume $L$ is geometric and $p\vee \sigma(p)=1$ for all elements $p$ of $L$. Since
\[r(p\vee \sigma(p)) + r(p\wedge \sigma(p)) \le r(p) + r(\sigma(p))\]
we see that $r(p)+r(\sigma(p))\ge d$ for all $p$. So if $r(p)\le k$ then $r(\sigma(p))\ge d-k$.\qed

Recently Huh and Katz \cite{HuhKatz} settled a longstanding open problem by proving that
the Whitney number $W_i$ for a geometric lattice form a 
log-concave (and hence unimodal) sequence.
As stated, Corollary~\ref{18.2} is due to Dowling and Wilson \cite{DowlingWilson}, but the most interesting case is when $k=1$, where the result was first established by Basterfield and Kelly \cite{BasterfieldKelly}, and independently by C. Greene \cite{Greene1970}. If equality holds in 
Corollary~\ref{18.2}, then Dowling and Wilson \cite{DowlingWilson} prove that $L$ must be modular; 
in the case $k=1$ this was also observed in \cite{BasterfieldKelly,Greene1970}. We take this further 
in Section~{sec:py-hyps}.

If $p$ is a point in $L$ then $\sigma(P)$ must be a complement to $p$. (This can be viewed as a consequence of the fact that points are modular elements.) Thus it is natural to ask if there could be a permutation $\sigma$ such that $\sigma(p)$ is a complement of $p$, for all elements $p$ in the lattice. This can be achieved under suitable conditions.

\begin{theorem}[Dowling \cite{Dowling1977}]\label{18.3}
	Let $L$ be a lattice such that $\mu_L(0,p)\mu_L(p,1)\ne 0$, for any element $p$. Then there is 
	a permutation $\sigma$ of $L$ such that $\sigma(p)$ is a complement of $p$, for all $p$ in $L$.
\end{theorem}

\proof
Let $G(p)$ denote the set of all elements $x$ of $L'$ such that $x\vee p<1$. Let $MN$ be the matrix with rows and columns indexed by $L$, such that 
\[(M)_{pq}:= \mu(G(p)_{\le q}).\]
By Lemma~\ref{lem:mucomp}, the $pq$-entry of $M$ is zero if $p$ and $q$ are not complements, so we can prove the theorem by showing that $\det M \ne 0$. (Perhaps it is worth noting that $M$ is probably not symmetric.)

If $x\in L'$ then $p\ne 0$ and so $\sum_{z\le p}\mu_L(0,z)=0$. Hence
\begin{align*}
0 &= \sum_{z\le q, z\vee p<1} \mu_L(0,z) + \sum_{z\le q, z\vee p=1} \mu_L(0,z)\\
&=\mu(G(p)_{\le q}) + \sum_{z\le q, z\vee p=1} \mu_L(0,z).
\end{align*}
If $H$ denotes the matrix we used in the proof of Theorem \ref{18.1} and $D$ is the diagonal matrix with $(D)_{pp}=\mu_L(0,p)$, then
\[M = -Z^T D H,\]
from which the theorem follows immediately.\qed

Unfortunately Theorem \ref{18.3} does not seem to lead to any strengthening of Corollary \ref{18.2}.

\section{Modular Lattices}

We have seen that if $L$ is geometric with height $n$ then $W_1\le W_{n-1}$. It is reasonable to ask what can be said if equality holds. As we will see in the next section, the answer is that $L$ must be modular. For this to make sense we must first define modular lattices themselves.

We start with an identity due to Dedekind, which holds in any lattice $L$. Suppose that $a,b$ and $c$ are elements of $L$. Then $a\vee (b\wedge c)$ lies below both $a\vee b$ and $a\vee c$. Hence $a\vee (b\wedge c)\le (a\vee b)\wedge (a\vee c)$ and so we see we have proved that, if $a\le c$ then
\begin{equation}\label{19.1}
a\vee (b\wedge c) \le (a\vee b)\wedge c.
\end{equation}

This is Dedekind's identity. A lattice is \textsl{modular} if equality holds in Dedekind's identity for all $a,b$ and $c$ in $L$ with $a\le b$. It is not hard to verify that any sublattice of a modular lattice is modular, and that products of modular lattices are modular. The dual of a modular lattice is modular. All modular lattice are ranked, but the proof of this is left to you as well.

\textbf{Exercise}: Show that a lattice is modular if and only if it is both upper and lower semimodular.

A word of warning here. We defined modular elements of geometric lattices. Note though that even if $a$ is a modular element of $L$, neither $[0,a]$ nor $[a,1]$ need be modular.

Our next lemma shows that we can test if a lattice is modular without looking at all triples $a,b$ and $c$ where $a\le c$.

\begin{lemma}
A lattice $L$ is modular if and only if $a\vee (b\wedge c)=c$ whenever $c\in [a, a\vee b]$.
\end{lemma}
\proof
Maybe this will be left for you too.
\qed
We now turn to characterization of modular geometric lattices.

\begin{lemma}
Let $L$ be a geometric lattice. An element $b$ of $L$ is modular if and only if $a\vee (b\wedge c) = (a\vee b)\wedge c$ whenever $a\le c$.
\end{lemma}
\proof
Suppose that $a\le c$ and let $w$ and $w'$ respectively denote $(a\vee b)\wedge c$ and $a\vee (b\wedge c)$. Note that $w'\ge w$, by Dedekind's identity. We have
\[b\wedge c \ge b\vee (a\vee (b\wedge c))\ge v\wedge (b\wedge c)=b\wedge c\]
and 
\[b\vee a \le b\vee ((a\vee b)\wedge c) \le b\vee (b\vee a),\]
therefore $v\wedge w'=b\wedge c$ and $b\vee w = b\vee a$. It is even easier to verify that $b\wedge w=b\wedge c$ and $b\vee w'=b\vee a$. So if $b$ is modular we have both
\[r(w) + r(b) + r(b\vee w) + r(b\wedge w) = r(b\vee a) + r(b\wedge c)\]
and 
\[r(w') + r(b) =r(b\vee w') + r(b\wedge w') + r(b\vee a) + r(b\wedge c),\]
whence $r(w) + r(w')$. As $w\le w'$, this implies that $w=w'$.

Now suppose that $a\vee (b\wedge c) = (a\vee b) \wedge c$ whenever $a\le c$. We show that no two complements of $b$ are comparable. But if $a$ and $c$ are complements to $b$ and $a\le b$ then
\[a = a\vee 0 = a\vee (b\wedge c) = (a\vee b)\wedge c = 1 \wedge c=c.\]
Hence the complements of $b$ form an antichain, and so $b$ is a modular elements.
\qed

\begin{corollary}
A geometric lattice is modular if and only if each element is modular.
\end{corollary}

\begin{lemma}
A geometric lattice is modular if and only if each hyperplane is a modular element of $L$.
\end{lemma}
\proof
Suppose $c\in L$ and $c\not\le h$. Then $h\vee c=1$ and $h\wedge c$ is a hyperplane in the interval $[0,c][$. We first show that if $h$ is modular in $L$ and $h\wedge c$ is modular in $[0,c]$. If $h$ is a hyperplane and $a\not\le h$ then
\[r(h) + r(a) \ge r(h \vee a) + r(h\wedge a) = r(a) + r(h\wedge a).\]
Hence $h$ is modular if and only if $r(a) - r(h\wedge a) = r(1) - r(h)=1$, i.e., if and only if $a\wedge h$ is covered by $a$, for any $a$ in $L$ such that $a\not\le h$. If $b\le a$ then $(h\wedge a)\wedge b = h\wedge b$. If $h$ is modular then $b$ must cover $h\wedge b$, whence $(h\wedge a)\wedge b$ is covered by $b$. If follows that $h\wedge a$ is a modular hyperplane in $[0,a]$.

Now suppose that $L$ is a geometric lattice in which every hyperplane is modular. We prove $L$ is modular by induction on its height. If $a<1$ then all hyperplanes in $[0,a]$ are modular. (You should show that if $a$ covers $nb$ then there is a hyperplane $h$ such that $h\wedge a=b$.) By induction it follows that all proper intervals of $L$ are modular.

Now let $a$ and $b$ by any two elements of $L$. We want to verify that 
\begin{equation}
r(a) + r(b) = r(a\vee v) + r(a\wedge b).
\label{19.2}
\end{equation}
If either $a\vee b>1 $ or $a\wedge b>0$ then this already holds, by our induction hypothesis. We therefore assume that $a$ and $b$ are complements. Since every element of $L$ is the intersection of the hyperplanes containing it, there is a hyperplane $h$ containing $a$ but not $b$. If $h\wedge b=0$ then, since $h$ is modular $b$ must be a point. As all points are modular, equality then holds in (\ref{19.2}).

Hence we may assume that $h\wedge b>0$. Since $[0,h]$ is modular 
\[r(a) + r(b\wedge h) = r(a\vee(b\vee h)) + r(a\wedge b \wedge h) = r(a\vee (b\wedge h))\]
and, since $[b\wedge h, 1]$ is modular
\[r(b)+r(a\vee(v\wedge h)) = r(1) + r(b\wedge (a\vee (v\wedge h))).\]
Combining these two inequalities we deduce that
\[r(b) + r(a) = r(1) + r(b\wedge (a\vee (b\wedge j))) - r(b\wedge h).\]
Let $w$ denote $b\wedge (a\vee (b\wedge h))$. Clearly $w\ge b\wedge h$, if we can show that $w=b\wedge h$ then $r(a)+r(b)=1$, as required.

Since $a\wedge b=0$, it is trivial to verify that $b\wedge h$ is a complement to $a$ in $[0,a\vee (b\wedge h)]$. Further
\[a\wedge w = a\wedge (b\wedge (a\vee()b\wedge h)))\le a\wedge b=0\]
while
\[a\vee (b\wedge h) =a\vee (a\vee (b\wedge h)) \ge a\vee (b\wedge (a\vee (b\wedge h)))\ge a\vee (b\wedge h)\]
and so $w$ is a second complement to $a$ in $[0,a\vee(b\wedge h)]$. Since this interval is contained in $[0,h]$ it is modular and, as $b\wedge h\ge w$ it follows that $w=b\wedge h$.
\qed

A \textsl{line} in a geometric lattice is an element of height two, i.e., the join of two points.

\begin{lemma}\label{19.5}
A hyperplane $h$ in a geometric lattice is modular if and only if $h\wedge \ell >0$, for every line $\ell$.
\end{lemma}
\proof
Let $h$ be a hyperplane. To show $h$ is modular we need only verify that $b\not\le h$ then $b$ covers $b\wedge h$. Assume by way of contradiction that $b\not\le h$ and $b$ does not cover $h$. We will use this to find a lines meeting $h$ in $0$.

Let $c$ be a complement to $b\wedge h$ in $[0,b]$. Since $r(b\wedge c)=0$ we have $r(c)\ge r(b)-r(b\wedge h)$, consequently $r(c)\ge 2$ and so there is a line $\ell$ lying below $c$. But then
\[h\wedge \ell = h\wedge (b\wedge \ell) = (h\wedge b)\wedge \ell =0.\]
It follows that $h$ is modular.
\qed

\section{Points and Hyperplanes (again)}
We will prove now that if the number of points in a geometric lattice $L$ is equal to the number of hyperplanes then $L$ is modular. Both proofs proceed by showing that if $L$ has rank $d$ and $W_1(L)=W_{d-1}(L)$ then any hyperplane meets any line non-trivially. (So there was some point to the trials of the previous section.) One reason this result is so interesting is that the structure of complemented modular lattice is very restricted: every complemented modular lattice is a direct sum of subspace lattices $\mathcal{B}_q(n)$ and copies of $\mathcal{B}(1)$.

\begin{theorem}[Basterfield and Kelly \cite{BasterfieldKelly}]\label{20.1}
Let $L$ be a geometric lattice of rank $d$. Then $W_1(L) = W_{d-1}(L)$ if and only if $a\vee b=1$. As we saw in dets, we have
\[G = ZFZ^T\]
where $Z = Z_L$ and $F$ is the diagonal matrix with $(F)_{aa}=\mu(a,1)$. Since $L$ is geometric, $F$ is invertible and therefore $G$ is invertible. We claim that if $a\wedge b=0$ then $(G^{-1})_{ab}\ne 0$. In fact we have $G^{-1}=(Z^T)^{-1} F^{-1} Z^{-1}$ and therefore
\[(G^{-1})_{ab} = \sum_{x\le a\wedge b} \frac{\mu(x,a)\mu(x,b)}{\mu(x,1)}.\]
When $a\wedge b=0$ this implies that
\[(G^{-1})_{ab} = \frac{\mu(0,a)\mu(0,b)}{\mu(0,1)}\ne 0.\]

We may write $G$ in partitioned form as
\[G = \pmat{0 & M\\N& X}\]
where the rows of $M$ are indexed by the points of $L$ and its columns by the hyperplanes. Since $G$ is invertible the rows of $M$ are linearly independent; this proves again that $W_1(L)\le W_{d-1}(L)$.

Now assume that $W_1(L) = W_{d-1}(L)$. Then $M$ and $N$ are square invertible matrices and accordingly
\[G^{-1} = \pmat{ -N^{-1}XM^{-1} & N^{-1}\\ M^{-1} & 0}.\]
What matters here is the zero submatrix of $G^{-1}$--its presence shows that if $a$ is not a point and $h$ is a hyperplane of $L$ then $G^{-1}_{ah}=0$. This implies that $h\wedge a>0$ and consequently for all hyperplanes $h$ and all lines $\ell$ of $L$ we have $h\wedge \ell >0$. Therefore $L$ is modular by Lemma \ref{19.5}.
\end{theorem}

The proof of Theorem \ref{20.1} can be extended to show that if equality holds in Corollary \ref{18.2} then $L$ is modular. (This is not an unreasonable exercise.) 

\begin{corollary}[Greene \cite{Greene1970}]
Let $L$ be a geometric lattice of rank $d$. Then $W_1(L)\le W_i(L)$, and if equality holds then $i=d-1$.
\end{corollary}
\proof
Let $L$ be a geometric lattice and let $f$ be the map from $L$ to $L$ defined by
\[f(x) =\begin{cases}
x,\quad \text{if } r(x)<k;\\
1,\quad \text{otherwise}.
\end{cases}\]
Then $f$ is order-preserving and its image is a geometric lattice (albeit, not a sublattice of $L$). The image of $L$ under $f$ is an \textsl{upper truncation} of $L$. Suppose that $W_1(L)=W_i(L)$ and let $L'$ be the geometric lattice obtained by truncating $L$ at height $i+1$. Then $L'$ is geometric and we have
\[W_1(L) = W_1(L') \le W_i(L') = W_i(L).\]
Further, if the first and last terms here are equal then $L'$ is modular, by the previous theorem.

Assume $i<d$ and choose an element $a$ of $L$ with rank $i-2$. Then the interval $[a,1]$ in $L$ has height four and therefore it contains a set of four independent points. These four points generate a sublattice isomorphic to $\mathcal{B})4)$, and the elements of rank two in it have rank $i$ in $L$. Hence the interval $[a,1]_{L'}$ has height three and contains four points and six lines. This implies that it is not modular, but every interval in a modular lattice is modular and therefore $L'$ cannot be modular.
\qed

Our proof of Theorem \ref{20.1} was based on the approach of Dowling and Wilson. We now present a version of the original proof of Basterfield and Kelly. (It is simple and elegant--our only criticism is that it does not use the M\"obius function.)

Assume that $L$ is a geometric lattice with rank $d$. We aim to prove by induction on $d$ that $W_1(L)\le W_{d-1}(L)$, with equality implying that $L$ is modular. Let $p$ be a point and $h$ a hyperplane of $L$ such that $p\not\le h$. We make two claims:
\begin{enumerate}[(a)]
\item $W_{d-2}[p,1]\ge W_1[0,h]$ and
\item $W_1[p,1]\ge W_1[0,h]$.
\end{enumerate}
To prove (a), suppose that $a$ and $b$ are covered by $h$. If $p\vee a = p\vee b$ then
\[p\vee a=p\vee a\vee b=p\vee h,\]
but since $p$ is modular $r(p\vee a)<r(p\vee h)$ and therefore the map $x\mapsto x\vee p$ is in injection from the hyperplanes of $[0,h]$ into the hyperplanes of $L$ on $p$. Thus $W_{d-2}[0,h]\le W_{d-2}[0,p]$. Since $[0,h]$ is geometric, $W_1[0,h]\le W_{d-2}[0,h]$ by induction and thus (a) is proved. For (b), the map $x\mapsto x\vee p $ is a bijection from the points of $[0,h]$ to the lines of $L$ on $p$ which intersect $h$ nontrivially.

Now we prove that $W_1(L)\le W_{d-1}(L)$ and that, if equality holds, $W_{d-2}[p,1]=W_1[0,h]$ for any point $p$ and hyperplane $h$ such that $p\wedge h=0$. Let $\comp{p}$ and $\comp{h}$ denote $W_{d-2}[p,1]$ and $W_1[0,h]$ respectively. By (a) above, $\comp{p}\ge \comp{h}$. Let $n$ be the number of points and $m$ be the number of hyperplanes in $L$. Then
\begin{equation}
n = \sum_p \frac{m-\comp{p}}{m-\comp{p}} = \sum_[p,h: p\wedge h=0] \frac{1}{m-\comp{p}}\ge \sum_{p,h: p\wedge h=0} \frac{1}{m-\comp{h}} = \sum_h \frac{n-\comp{h}}{m-\comp{h}}.
\label{eqn20.1}
\end{equation}
If $m<n$ then
\[\frac{n-\comp{h}}{m-\comp{h}}> \frac{n}{m}\]
whence the last term in (\ref{eqn20.1}) is strictly greater than $n$. Hence we conclude that $m\ge n$ and, if equality holds, $\comp{p}=\comp{h}$ for any point $p$ and hyperplane $h$ with $p\wedge h=0$.

Assume that $m=n$ and that $p$ is a point and $h$ is a hyperplane not on $p$. Since $[p,1]$ is geometric, we may use (b) above to deduce that 
\[\comp{h}=W_1[0,h]\le W_1[p,1]\le W_{d-2}[1,p]=\comp{p}.\]
This implies that $W_1[0,h] = W_1[p,1]$, and our proof of (b) then implies that every line on $p$ meets $h$ non-trivially. Thus we have shown that if $h$ is a hyperplane and $\ell$ is a line in $L$ then $h\wedge \ell >0$ and therefore $L$ is modular, by Lemma \ref{19.5}.

\section{Kung Fu?}
We will describe some important work of J. Kung \cite{Kung1985}, in a formulation communicated privately to the author by C. Greene. This provides yet another approach to some of the work in hypers and pt-hyps.

If $f$ is a function on a lattice $L$, let $\hat{f}$ be define by 
\[\hat{f}(a) = \sum_{x\le a} f(x).\]
Our main theorem can be viewed as providing one answer to the following problem. Suppose that $A$ and $B$ are subsets of a lattice $L$. What conditions on $A$ and $B$ guarantee that any function $f$ on $L$ with support in $A$ is determined by the restriction $\hat{f}\upharpoonright B$ of $\hat{f}$ to $B$? (Admittedly this appears to be a convoluted problem, with little hope of a useful answer arising.)

\begin{theorem}[Kung \cite{Kung1985}]
Let $A$ and $B$ be subsets of the lattice $L$ such that, if $x\in :$ then either
\begin{enumerate}[(a)]
\item 
$x\in B$, or
\item 
there exists $x^*$ in $L$ such that $\mu(x,x^*)\ne 0$ and $a\vee x\ne x^*$ if $a\in A$. Then $\hat{f}\upharpoonright B$ determines $\hat{f}$, and there is an injection $\phi: A\mapsto B$ such that $\phi(a)\ge a$ for all $a$ in $A$.
\end{enumerate}
\end{theorem}

Before embarking on the proof of this result, we present one application. Let $L$ be a geometric lattice with rank $d$, let $A$ be the set of elements of $L$ with rank at most $k$ and let $B$ be the set of elements with rank at least $d-k$. If $x\in L$, let $x^*=1$.

If $x\notin B$ then $\mu_L(x,1)\ne 0$. If, further, $a\in A$ then
\[r(x\vee a) \le r(x) + r(a) \le d-k-1+k <d=r(1).\]
Hence the conditions of Theorem \ref{20.1} are (and be be!) satisfies. What may we conclude? Let $f_a$ be the function on $L$ defined by
\[f_a(x)=\begin{cases}
1,\quad \text{if } x=1;\\
0,\quad \text{otherwise}.
\end{cases}\]
Then $\hat{f}_a(b) =1$ if $b\ge a$, and is zero otherwise. The theorem implies that, for functions $f$ supported on $A$, the linear mapping
\[f\mapsto \hat{f}\upharpoonright B\]
is injective. This implies that $\dim\re^B\ge \dim \re^A$, from which it follows that $\abs{A}\le \abs{B}$. This provides another proof of Corollary \ref{18.2}. In fact a stronger statement can be made. The function $\hat{f}_a$ can be identified with the row of $Z_L$ indexed by $a$ and therefore Kung's theorem implies that the submatrix of $Z_L$ with rows indexed by elements of $A$ and columns by elements of $B$ has linearly independent rows. Hence there is an injection $\phi:A\to B$ such that $\phi(a)\ge a$, for all $a$ in $A$.

We start the proof of Theorem \ref{20.1} now. If $x$ and $y$ are elements of $L$ such that $x<y$ 
and $f$ is a function on $L$, we have
\begin{align*}
	\sum_{t\in [x,y]} \mu_L(t,y)\hat{f}(t) 
		&=\sum_{t\in [x,y]}\sum_{s\le t} \mu_L(t,y)f(s) \\
		&= \sum_s f(s)\sum_{t\in [x\vee s,y]} \mu_L(t,y)\\
		&=\sum_{s: s\vee x=y} f(s).
\end{align*}
Now suppose $x\notin B$ and let $y$ be $x^*$. Then, if the support of $f$ is contained in $A$, 
the last term above is zero and so
\begin{equation}
	\mu_L(x, x^*)\hat{f}(x) = -\sum_{x<t \le x^*} \mu_L(t,x^*) \hat{f}(t).
\label{eq:muL}
\end{equation}
Condition (b) of the theorem implies that $1$ must lie in $B$. By Equation~\ref{eq:mUL}), if $x\notin B$ 
then $\hat{f}(x)$ is determined by the value of $\hat{f}$ on elements $t$ in $L$ such that $t>x$. 
It follows by induction that $\hat{f}$ is determined by $\hat{f}\upharpoonright B$. This completes 
the proof.

We describe a second application of Theorem~\ref{21.1}, to modular lattices. An element $a$ in a lattice $L$ is \textsl{join-irreducible} if, whenever $x\vee y=a$, either $x=a$ or $y=a$. In other words, $a$ is not the join of two smaller elements. The set of all join-irreducible elements of $L$ will bve denoted by $J(L)$. Similarly we may define \textsl{meet-irreducible} elements; the set of meet-irreducible elements of $L$ will be denoted by $M(L)$. Note $0\in J(L)$ and $1\in M(L)$, hence these subsets are not empty. In a geometric lattice, $J(L)$ consists of $0$ and all the points, while $M(L)$ consists of the hyperplanes and $1$. 

Assume $L$ is modular with $A=J(L)$ and $B=M(L)$. If $x\notin B$, define $x^*$ to be the join of the elements which cover $x$. Then $[x,x^*]$ is a modular point-lattice, therefore it is geometric and $\mu_L(x,x^*)\ne 0$. Suppose $a$ is join-irreducible. Since $L$ is modular, the intervals $[a\wedge x,a]$ and $[x,x\vee a]$ are isomorphic (by Lemma \ref{13.1}). But this implies that $x\vee a$ is join-irreducible in $[x,x\vee a]$, which is a geometric lattice. Therefore $x\vee a$ must cover $x$, and hence cannot be equal to $x^*$. Thus the conditions of Kung's theorem are satisfied, whence we conclude that in a modular lattice $\abs{J(L)}\le \abs{M(L)}$. As $L^{\op}$ is modular if $L$ is and 
\[
	J(L^{\op})=M(L),\quad M(L^{\op})=J(L)
\]
it follows that $\abs{J(L)}=\abs{M(L)}$ for modular lattices. (This is a famous result of Dilworth's.)

\section{Contraction and Deletion}

Let $L$ be a point lattice with point set $\Omega$ and suppose $p\in \Omega$. We define a function $f$ from $L$ into itself as follows:
\[
	f(a) = \vee \{q:q\in \Omega\backslash p, q\le a\},
\]
with the understanding that $f(0)=0$. It is easy to check that $f$ is order preserving and that $f(x)$ is either $x$ itself or the unique element covered by $x$ and not in $[p,1]$. Hence $f$ is a decreasing map. Note that $f$ is an order preserving and decreasing map from $L'\backslash p$ into itself.

Let $M$ be the poset of fixed points of $f$. This is a join semi-lattice with a $0$- and $1$-element. (The latter is usually the $1$-element of $L$.) If $M':=M\backslash \{0,1_L\}$ then, by 
Lemma~\ref{lem:avmu0}, $M'$ is a retract of $L'\backslash p$ and therefore $\mu(M')=\mu(L'\backslash p)$. On the other hand by Lemma \ref{14.1} we have
\[\mu(L') = \mu(L'\backslash p) + \mu(L'_{<p})\mu(L'_{>p}).\]
Since $p$ is a point, $\mu(L'_{<p})=-1$ and therefore
\[\mu_L(0,1) = \mu(M') - \mu_L(p,1).\]
There are two cases to be considered. If the join of the points of $L$ distinct from $p$ is equal to $1_L$ then $\mu(M')=\mu_M(0,1)$. If the join of the points of $L$ distinct from $p$ is not equal to $1_L$ then, because of the careful way we defined it, $M'$ has a $1$-element and $\mu(M')=0$. We will call a point $p$ a \textsl{co-loop} if the join $h$ of the points distinct from $p$ is not equal to $1$. (If $L$ is geometric and $h\ne 1$ then $h$ is a modular hyperplane.) Since $M$ is a semi-lattice it gives rise naturally to a lattice that we will denote by $L\backslash p$. (This is not a particularly good choice of notation, but will do for now.) We can summarise our conclusions as follows.

\begin{lemma}\label{22.1}
Let $L$ be a point lattice and let $p$ be a point of $L$. Then
\[\mu_L(0,1) = \begin{cases}
-\mu_L(p,1),\quad \text{if $p$ is a co-loop};\\
\mu_{L\backslash p}(0,1)-\mu_L(p,1),\quad \text{otherwise}.
\end{cases}\]
\end{lemma}
If $L$ is the lattice of contractions of a graph $G$ and $e$ is an edge in $G$ then $e$ is a point in $L$ and $L\backslash e$ is the lattice of contractions of the graph $G\backslash e$, obtained by deleting $e$ from $G$. The interval $[e,1]$ in $L$ is the lattice of contractions of $G'/e$, which is the graph obtained from $G$ by contracting the edge $e$. If can be shown that if $L$ is geometric then so is $L\backslash p$.

\textbf{Exercise}: Use Lemma \ref{22.1} to prove the broken circuit theorem (Theorem \ref{15.1}).

The only significant application of Lemma \ref{22.1} I know of is to geometric lattices. There are many other classes of point lattices though---the face lattices of convex polytopes, for example.

\section{Null Designs}
Let $P$ be a poset. A function of \textsl{strength} at least $t$ on $P$ is a function $f$ with values in some ring such that 
\[\sum_{x\ge a} f(x) =0.\]
for any element $a$ of $P$ with height at most $t$. (In practice we assume that the ring is the ring of integers.) The most important case is when $P$ is the lattice of all subsets of a set $V$, when a function of strength at least $t$ is sometimes called a null $t$-design. Our basic problem is to find good lower bounds on the support of a function of strength $t$.

If $f$ is any function on $P$ and $\hat{f}$ is as in the previous section then $f$ has strength at least $t$ if and only if the support of $\hat{f}$ contains no elements of height $t$ or less. Let $P$ be a poset with zero. If $b\in P$ let $f_b$ be the function obtained by M\"obius inversion on $[0,b]$ to $\hat{f}$. Then $f_b$ is a function on $[0,b]$, which we extend to a function on $P$ by setting $f_b(x)=0$ if $x\not\le b$. It is immediate that $f_b$ is a function of strength at least $t$ with support contained in $[0,b]$. We have two formulas for computing the values of $f_b$.

\begin{lemma}
Let $P$ be a meet semi-lattice. If $b\in P$ then
\[f_b(c)=\sum_{x\wedge b=c} f(x).\]
\end{lemma}
\proof
We have
\begin{align*}
f_b(c)
&=\sum_{y\le b} \mu(c,y)\hat{f}(y)\\
&=\sum_{y\le b} \mu(c,y) \sum_{x\ge y} f(x)\\
&=\sum_{x,y : c\le y\le x} \mu(c,y)f(x)\\
&=\sum_x \left(\sum_{y:c\le y\le x} \mu(c,y)\right) f(x)\\
&=\sum_{x\wedge b=c}f(x).
\end{align*}
The next result is trivial to verify.

\begin{lemma}\label{23.2}
Let $f$ be a function on the poset $P$ with strength at least $t$ and let $b$ be an element of $P$ in the support of $\hat{f}$ with minimal height. If $c\le b$ then
\[f_b(c) =\mu(c,b)\hat{f}(b).\tag*{\sqr53}\]

As an immediate corollary of Lemma \ref{23.2}, we see that the support of $f$ is bounded below by
\[\abs{\{c\le b: \mu(c,b)\ne 0\}}.\]
This bound can be improved in two cases. The previous two lemmas combine to yield that
\begin{equation}\label{eqn23.1}
\mu(c,b)\hat{f}(b)=\sum_{x\wedge b=c}f(x).
\end{equation}
\end{lemma}

\begin{lemma}
Let $P$ be a meet semi-lattice and let $f$ be a $(0,\pm 1)$-valued function of strength $t$ on $P$. If $b$ is an element of height $t$ in $P$ such that $\hat{f}(b)\ne 0$ then the support of $f$ has size at least
\[\sum_{c\le b}\abs{\mu(c,b)}.\]
\end{lemma}
\proof
If $f$ is $(0,\pm 1)$-valued then $\abs{\hat{f}(b)}\ge 1$ and (\ref{eqn23.1}) implies that there are at least $\abs{\mu(c,b)}$ elements $x$ in $P$ such that $x\wedge b=c$ and $f(x)\ne 0$.\qed

Both the previous lemma and the following theorem seem to have appeared first in unpublished work of Cho.

\begin{theorem}
Let $P$ be a semi-lattice, let $f$ be a function on $P$ with strength $t$ which is supported by elements of height $t+1$. If $b$ is an element of $p$ with height $t+1$ such that $\hat{f}\ne 0$ then the support of $f$ has size at least
\[\sum_{c\le b} \abs{\mu(c,b)}.\]
If equality holds then $f$ is $(0,\pm 1)$-valued.
\end{theorem}

\proof
Assume that, amongst all elements of height $t+1$ in the support of $\hat{f}$, we have chosen $b$ so that $\hat{f}(b)=f(b)$. By (\ref{eqn23.1}),
\[\mu(c,b) = \sum_{x\wedge b=c} \frac{f(x)}{f(b)}\le \sum_{x:x\wedge b=c,f(x)\ne 0} 1.\]
The theorem follows at once.\qed

It is reasonable to ask which meet semi-lattices we would like to apply the results of this section to.
The first two cases of interest are $\mathcal{B}(n)$ and $\mathcal{B}_q(n)$, which are lattices. We also
have the subspace lattice of a polar space. Finally suppose that $V$ is a $d$-dimensional vector space
over a finite field and $U$ is a subspace of $V$. Then the set of subspaces of $V$ which intersect $U$
in the zero subspace is a meet semi-lattice. Both these last examples have the property that any
interval is the subspace lattice of a projective space, and thus its M\"obius function is known.

For $\mathcal{B}(n)$, it is not too hard to prove that a function of strength at least $t$ has support of size at least $2^{t+1}$. This is stronger than we have just proved. But for $\mathcal{B}_q(n)$ the results of this section are the strongest known. It is somtimes possible to give a simpler expression for
\[\sum_{c\le b} \abs{\mu(c,b)}.\]
If $P=\mathcal{B}(n)$ and $b$ is a set of size $(t+1)$, this sum equals $2^{t+1}$. If $P$ is $\mathcal{B}_q(n)$ and $b$ has dimension $t+1$ then it equals
\[\prod_{i=0}^t (1+q^i).\]
If $P=\mathcal{P}(n)$ and $b$ is a partition with exactly $n-k$ cells then our sum equals $(n-k)!$. (These claims all follow from \cite[Proposition 4.20]{Aigner1996}.)

\bibliographystyle{amsplain}

\end{document}